\documentclass{amsart}
\usepackage{amsfonts}
\usepackage{amscd}
\usepackage{amsmath}

\usepackage[all]{xy}
\CompileMatrices

\theoremstyle{plain} \numberwithin{equation}{section}
\newtheorem{theorem}{Theorem}[section]
\newtheorem{corollary}[theorem]{Corollary}
\newtheorem{conjecture}{Conjecture}
\newtheorem{lemma}[theorem]{Lemma}
\newtheorem{proposition}[theorem]{Proposition}
\theoremstyle{definition}
\newtheorem{definition}[theorem]{Definition}

\newtheorem{remark}[theorem]{Remark}

\numberwithin{equation}{section}

\newcommand{\cO}{\mathcal{O}}
\newcommand{\cH}{\mathcal{H}}
\newcommand{\cK}{\mathcal{K}}
\newcommand{\cT}{\mathcal{T}}
\newcommand{\cU}{\mathcal{U}}
\newcommand{\cW}{\mathcal{W}}

\newcommand{\ad}{\mathrm{ad}}

\newcommand{\fa}{\mathfrak{a}}

\newcommand{\fc}{\mathfrak{c}}
\newcommand{\fg}{\mathfrak{g}}
\newcommand{\fh}{\mathfrak{h}}
\newcommand{\fk}{\mathfrak{k}}

\newcommand{\fm}{\mathfrak{m}}
\newcommand{\fn}{\mathfrak{n}}
\newcommand{\fp}{\mathfrak{p}}
\newcommand{\fq}{\mathfrak{q}}
\newcommand{\fs}{\mathfrak{s}}
\newcommand{\ft}{\mathfrak{t}}

\newcommand{\fz}{\mathfrak{z}}
\newcommand{\ga}{\alpha}

\newcommand{\gD}{\Delta}
\newcommand{\gO}{\Omega}
\newcommand{\cL}{\mathcal {L}}

\newcommand{\R}{\mathbb{R}}
\newcommand{\C}{\mathbb{C}}

\newcommand{\Z}{\mathbb{Z}}

\newcommand{\Hsp}{\cH^2(\Xi_+)}
\newcommand{\pr}{\mathrm{pr}}

\newcommand{\SL}{\mathrm{SL}}
\newcommand{\Ad}{\mathop{\mathrm{Ad}} }

\newcommand{\id}{\mathrm{id}}
\newcommand{\ev}{\mathrm{ev}}

\newcommand{\im}{\mathop{\rm im}}

\newcommand{\oline}{\overline}
\newcommand{\la}{\langle}
\newcommand{\ra}{\rangle}

\begin{document}
\title[Horospherical model for discrete series]{Horospherical model for
holomorphic discrete series and horospherical Cauchy transform}

\author{Simon Gindikin, Bernhard Kr\"otz, and Gestur {\'O}lafsson}\thanks{The
first author was supported by the Louisiana Board of
Regents grant \textit{Visiting Experts in
Mathematics}. The last two authors were
supported by the Research in Pairs program of the Mathematisches
Forschungsinstitut, Oberwolfach. The last author
was supported by NSF grant DMS-0139783 and DMS-0402068}
\address{Department of Mathematics, Rutgers University, New Brunswick, NJ 08903, USA}
\email{gindikin@math.rutgers.edu}
\address{Research Institute for Mathematical Sciences,
Kyoto University, Kyoto 606-8502, Japan}
\email{kroetz@kurims.kyoto-u.ac.jp}
\address{Department of Mathematics, Louisiana State
University, Baton Rouge,
LA 70803, USA} \email{olafsson@math.lsu.edu}

\maketitle

\section*{Introduction}

\noindent
For some homogeneous spaces the method of horospheres delivers an
effective way to decompose representations in irreducible ones.
For Riemannian symmetric spaces $Y=G/K$
horospheres are
orbits of maximal unipotent subgroups of
$G$. They are parameterized by points of the horospherical
homogeneous space $\Xi_\R= G/MN$ where $N$ is a fixed  maximal
unipotent subgroup and $M=Z_K(A)$ as usual.
The horospherical transform maps sufficiently regular
functions on $Y$ to the corresponding average along the
horospheres. The crucial point is, that the abelian group
$A$ acts on $\Xi$ and that this action
commutes with the action of $G$.  The decomposition
of the natural representation of $G$ in $L^2(\Xi)$ in irreducible ones reduces to
the decomposition relative to $A$. In this way we obtain
all unitary spherical representations on $Y$ (with constant
multiplicity), except the complementary series.
The computation of the Plancherel measure on $Y$ is
equivalent to the inversion of the horospherical transform.

The method of horospheres works for several other types of  homogeneous
spaces, including complex semisimple Lie groups (considered as
symmetric spaces) but it has very serious restrictions: discrete
series representations lie in the kernel of the horospherical
transform, as well as all representations induced from
parabolic subgroups that are not minimal. In short, the
kernel is the orthocomplement of the most continuous
part of the spectrum. The simplest example when the
horospherical transform can not be inverted is for the group
$\SL(2,\R)$. In \cite {g00,g02,g04} a modification
of the method of horospheres was suggested: complex horospherical transform
(the horospherical Cauchy-Radon transform). For a homogeneous
space $X$ we consider the complexification $X_\C$ and instead
of real horospheres on $X$ we consider complex horospheres on
$X_\C$ without real points (they do not intersect $X$). The
integration along a real horosphere is equivalent to the
integration of a $\delta$-function on $X$ with support on
this horosphere. In the complex version we replace this
$\delta$-function by a Cauchy type kernel with singularities
on the complex horosphere without real point. In \cite {g00,g02} it
is shown that such a complex horospherical transform has no kernel for
$\SL(2;\R)$ and that it reproduces the Plancherel formula; in \cite {g04}
it is shown for all compact symmetric spaces.

\par The objective of this paper is to show that the complex horospherical
transform has no kernel on the holomorphic discrete series.
Holomorphic discrete series exist for affine symmetric spaces
$X=G/H$ of Hermitian type; $G$ is here a group of Hermitian type
\cite {hc55,oo91} . The corresponding part of $L^2(X)$ can be realized as
boundary values of Hardy space $\cH^2(D_+)$  in a Stein tube
$D_+\subset X_\C$ with edge $X$
\cite{hoo91}. Our aim is
to  define a
complex horospherical transform which has no kernel on holomorphic
$H$-spherical representations.

\par The first step is a construction of the space
that is going to be the image of the complex
horospherical transform. For that, we consider those complex horospheres in the
Stein symmetric space $X_\C=G_\C/H_\C$ that
are parameterized by points of the complex horospherical space
$\Xi=G_\C/M_\C N_\C$. In $\Xi$ we then consider an orbit $\Xi_+$ of
$G\times \cT_+$ where $\cT_+$ is an abelian semigroup in the complex
torus $T_\C=AT$ with the compact torus
$T$ as the edge. The space  $\cO(\Xi_+)$ of holomorphic functions
on $\Xi_+$
is the Fr\'echet model of the holomorphic discrete series. More exactly,
if we decompose this representation  with respect
to the compact torus $T$ we obtain $G$-modules which are lowest weight modules
(if they are irreducible); we obtain all such modules with
multiplicity one.
Using the abelian semigroup $\cT_+$ we can define a
Hardy type space $\cH^2(\Xi_+)$ with spectrum
"almost all" of the holomorphic discrete series.

\par The next step is a geometrical background for the construction of
the horospherical transform. Firstly, we prove that the
horospheres $E(\xi)$ parameterized by points $\xi \in \Xi_+$ do
not intersect $X$. We construct a simple Cauchy type kernel which
has no singularities on $X$ and the edge of its singularities
coincides with $E(\xi)$. Using this kernel we define the
horospherical Cauchy transform from $L^1(X)$ to $\cO(\Xi_+)$  which
can be extended on $L^2(X)$. The horospherical
transform decomposed under
$T$ yields the holomorphic spherical Fourier transform.

\par The last step is the inversion of the horospherical Cauchy
transform. We give the Radon type inversion formula using results
from \cite{k} for the holomorphic discrete series. Let
us remark that for $X=\SL(2,\R)$ the inversion formula was
obtained  in \cite{g00,g02} with tools from integral geometry on
quadrics. This method automatically extends on any symmetric
spaces of Hermitian type of
rank 1, i.e., the hyperboloids of signature $(2,n)$. Let us also pay
attention to the complete similarity of formulas of this paper and
formulas in \cite {g04} for compact symmetric spaces. It confirms
the view that finite-dimensional spherical representations are
similar to representations of holomorphic discrete series.

\section{Symmetric spaces of Hermitian type}\label{s-one}

\noindent
The objective of this section is to set up a standard
choice of terminology that will be used throughout the text.

\medskip  Let us fix
some  conventions upfront. For a real Lie algebra $\fg$
let us denote by  $\fg_\C=\fg\otimes_\R \C$ its
complexification. Likewise, if not stated otherwise,
for a connected Lie group $G$ we write $G_\C$
for its universal complexification.
If $\varphi: G\to H$ is a homomorphism of connected Lie groups, then we
will also denote by $\varphi$
\begin{itemize}
\item the derived homomorphism $d\varphi({\bf 1}): {\rm Lie}(G)
\to {\rm Lie}(H)$,
\item the extension of $\varphi$ to a holomorphic homomorphism
$G_\C\to H_\C$.
\end{itemize}

\smallskip Let $G$ be a connected semisimple Lie group with Lie algebra
$\fg$. We assume that $G\subset G_\C$ and that
$G_\C$ is simply connected.
Let $\tau :G\to G$ be a non-trivial involution and write
$H$, resp. $H_\C$, for the $\tau$-fixed points in $G$, resp. $G_\C$.
The object of concern is the affine symmetric space
$X=G/H$. We observe that $X$ is contained in
its complexification $X_\C =G_\C/H_\C$ as a totally real
submanifold. Write $x_0=H_\C$ for the base point in $X_\C$.

\par Let $\fh$ be the Lie algebra of $H$ and note that
$\fg=\fh+\fq$ with $\tau|_{\fq}= -\id_{\fq}$.
The symmetric pair $(\fg,\fh)$ is called
irreducible if $\fg$ does not contain any
$\tau$-invariant ideals except the trivial ones,
$\{0\}$ and $\fg$. In that case, either $\fg$ is
simple or $\fg=\fg_1\times \fg_1$ with $\fg_1$ simple
and $\tau (x,x')=(x',x)$. We say that $X$ is
irreducible, if $(\fg,\fh)$ is irreducible. From now on
we will assume, that $X$ is irreducible.

\par Fix a Cartan involution
$\theta :G\to G$ commuting with $\tau$. Denote by
$K<G$ the subgroup of $\theta$-fixed points  and
write $Y=G/K$ for the associated
Riemannian symmetric space.
Write $\fk$ for the Lie algebra of $K$.
Then $\fg=\fk+\fs$ with $\theta|_{\fs}=-\id_{\fs}$.
Notice that the universal complexification
$K_\C$ of $K$ naturally identifies with the $\theta$-fixed points
in $G_\C$.

\par We will assume that $G$ is a Lie group
of Hermitian type, i.e. $Y$ is Riemannian symmetric space
of Hermitian type. The assumption can be phrased
algebraically: $\fz(\fk)\neq \{0\}$ with $\fz(\fk)$ the center of
$\fk$.

\par We assume that $\tau$  induces an  anti-holomorphic involution
on $Y$ and then call $X$ an {\it affine symmetric space
of Hermitian type}.

\begin{remark} (a)  Our assumptions on $G$ and $\tau$ can be
phrased algebraically, namely:
$$\fz(\fk)\cap\fq\neq\{0\}\ .\leqno{\rm (A)}$$
Let us mention that
another way to formulate (A) is to say that $\fq$ admits
an $H$-invariant regular elliptic cone, i.e.
$X$ is compactly causal \cite{ho}.
\par \noindent (b) Symmetric spaces of Hermitian type resemble
compact symmetric spaces on an analytical level. Combined they form
the class of symmetric spaces which admit lowest weight modules
in their $L^2$-spectrum (holomorphic discrete series).
\end{remark}

\medskip Since $X$ is irreducible, it follows that
$\fz(\fk)\cap \fq=i\R Z_0$ is one dimensional.
It is possible to normalize $Z_0$ in such a way that
the spectrum of $\ad (Z_0)$ is
$\{-1, 0,1\}$. The zero-eigenspace is  $\fk_\C$.
We denote the $+1$-eigenspace in $\fs_\C$ by $\fs^+$, and the $-1$-eigenspace
by $\fs^-$.

\medskip Let $\ft$ be a maximal abelian subspace
in $\fq$ containing $iZ_0$. Then $\ft$ is contained in $\fk\cap \fq$.
Set $\fa=i\ft$ and note that $\fa_\C=\ft_\C$.

\par Let $\gD$ be the set of roots of $\ft_\C$ in $\fg_\C$,
$$\gD_n=\{\ga\in\gD\mid \fg_\C^\alpha \subseteq \fs_\C \}
=\left\{\ga\in \gD\mid \ga (Z_0)\in \{-1, 1\}\right\}$$
and
$$\gD_k=\{\ga\in\gD\mid \fg_\C^\alpha\subseteq \fk_\C\}
=\{\ga\in\gD\mid \alpha (Z_0)=0\}\, .$$
Then $\gD=\gD_k\dot{\cup}\gD_n$. The elements of $\gD_n$ are
called \textit{non-compact roots}, and the elements in
$\gD_k$ are called \textit{compact roots}. We choose an ordering
in $i\ft^*$
such that $\ga (Z_0)>0$ implies that $\ga\in \gD_n^+\subseteq \gD^+$.
Let ${\mathcal W}$ be the Weyl group of $\gD$ and ${\mathcal W}_k$ the
subgroup generated by the reflections coming from the compact roots.
As $s(Z_0)=Z_0$ for
all $s\in {\mathcal W}_k$, it follows that $\gD^+_n$ is ${\mathcal W}_k$-invariant.

\subsection{Polyhedrons, cones and the minimal tubes}\label{ss=11}
Set $A=\exp(\fa)$, $A_\C=\exp(\fa_\C)$,
$T=\exp (\ft)$ and $T_\C=\exp (\ft_\C)$.
We note that
$$A_\C=T_\C= TA\simeq T\times A\, .$$

\par For $\alpha\in \Delta$ let $\check\alpha\in \fa$ be its coroot, i.e.
$\check\alpha\in [\fg_\C^\alpha, \fg_\C^{-\alpha}]\cap \fa$
and $\alpha(\check\alpha)=2$.
Then
\begin{equation}\label{eq=om}
\Omega=\sum_{\alpha\in \Delta_n^+} \R_{>0} \cdot  \check\alpha
\end{equation}
defines a ${\mathcal W}_k$-invariant open convex  cone in $\fa=i\ft$ which contains $Z_0$. Often one refers
to $\Omega$ as the {\it minimal cone} (it is denoted $c_{\rm min}$
in \cite{ho}). Let us
remark that one can characterize $\Omega$ as
the smallest $\cW_k$-invariant open convex cone in $\fa$ which contains
a long non-compact coroot, i.e.
\begin{equation}\label{eq=br} \Omega= {\rm co}
\left({\mathcal W}_k (\R_{>0}\cdot\check\alpha)\right)\qquad (\alpha
\ \hbox{long in}\ \Delta_n^+)\, .
\end{equation}
 Here ${\rm co}(\cdot)$ denotes the convex hull of
$(\cdot)$.
\par
We set $A_+=\exp(\Omega)$ and note that
$A_+\subset A$ is an open semigroup.
Moreover
$$\cT_+=T\exp (\Omega)=TA_+\subset T_\C\, $$
defines a semigroup and complex polyhedron with edge $T$.
We also use the notation $A_-=\exp(-\Omega)$ and
$\cT_-= TA_-$.

\par Define  $G$-invariant subsets of
$X_\C$ by
$$D_{\pm}=G A_\pm\cdot x_0\subset X_\C\, .$$
According to \cite{n99} $D_+$ and $D_-$  are Stein domains in $X_\C$
with $X=G\cdot x_0$ as Shilov boundary. Subsequently we will refer
to $D_+$ and $D_-$  as {\it minimal tube in $X_\C$ with edge $X$}.

\subsection{Minimal $\oline\theta\tau$-stable parabolics}
Denote by $g\mapsto \oline g$ the complex conjugation
in $G_\C$ with respect to the real form $G$.
Let
$$\fn^+_\C=\bigoplus_{\ga\in\gD^+_k}\fk_\C^\ga\quad\mathrm{and}
\quad \fn^-_\C=\bigoplus_{\ga\in\gD^+_k}\fk_\C^{-\ga}\, .$$
Set
$$\fn_\C=\fn_\C^+\oplus \bigoplus_{\ga\in\gD^+_n}\fg_\C^\ga=\fn_\C^+\oplus \fs^+\, ,$$
$$\fm_\C=\{U\in \fh_\C\mid  (\forall V\in\ft) \  [U,V]=0\}\, ,$$
and
$$\fp_\C= \fm_\C\oplus \ft_\C \oplus \fn_\C \, .$$
Notice, that $\fm_\C$ is contained in $\fk_\C$, as $Z_0\in \ft_\C$.
The Lie algebra $\fp_\C$ is a \textit{minimal} $\oline \theta\tau$\textit{-stable
parabolic subalgebra of} $\fg_\C$.
Define subgroups of $G_\C$ by $M_\C=Z_{H_\C}(\ft_\C)\subset K_\C$,
and $N_\C=\exp (\fn_\C)$.

Note that $T_\C=A_\C$. Then the  prescription
$$P_\C =M_\C A_\C N_\C= M_\C T_\C N_\C$$
defines a
\textit{minimal} $\oline\theta\tau$\textit{-stable parabolic subgroup of} $G_\C$
whose Lie algebra is $\fp_\C$. Write $\Gamma=M_\C \cap A_\C =M\cap T$
and observe that $\Gamma$ is a finite $2$-group.
The isomorphic map
$$(M_\C \times _\Gamma A_\C)\times N_\C \to P_\C, \ \
([m,a],n)\mapsto man$$
yields the structural decomposition of $P_\C$.

\smallskip We denote by  $\ft\subseteq \fc$  a
$\tau$-stable Cartan subalgebra of $\fg$
contained in  $\fk$.
Then $\fc = \ft\oplus \fc_h$, where $\fc_h=\fc \cap \fh$.
Denote by $\Sigma$ the set of roots of
$\fc_\C$ in $\fg_\C$. Similarly we set
$\Sigma_n$, the set of non-compact roots, $\Sigma_k$, the
set of compact roots. We choose a positive system $\Sigma^+$ such that
$\Sigma^+|_{\ft} \backslash \{0\}=\Delta^+$.
\par Define tori in $G$ by $C=\exp \fc$ and $C_h=\exp\fc_h$.
We note that $C=TC_h\simeq T\times_\Gamma C_h$.

\section{Complex Horospheres I: Definition and basic properties}\label
{section=hor-I}

The objective of this section is to discuss (generic) horospheres
on the complex symmetric space $X_\C=G_\C/H_\C$. We will
show that the space of horospheres is $G_\C$-isomorphic
to the homogeneous space $\Xi=G_\C / M_\C N_\C$.
Further we will introduce a $G$-invariant subdomain $\Xi_+\subset \Xi$
which will be a central object for the rest of this paper.

\par Set
$$\Xi= G_\C/ M_\C N_\C\  $$
and write  $\xi_0=M_\C N_\C$ for the base point of $\Xi$.
Usually we express elements
$\xi\in \Xi$ as $\xi=g\cdot \xi_0$
for $g\in G_\C$.
\par Consider the double fibration
\begin{equation}\label{eq=df}
 \xymatrix { & G_\C/ M_\C \ar[dl]_{\pi_1} \ar[dr]^{\pi_2} &\\
\Xi & & X_\C\,.}\end{equation}
By a {\it horosphere} in $X_\C$
we understand a subset of the
form
\begin{equation}\label{eq=E} E(\xi)=\pi_2(\pi_1^{-1}(\xi)) \qquad (\xi\in \Xi)\, .
\end{equation}
For $\xi=g\cdot \xi_0$ we record that
$$E(\xi)=gM_\C N_\C\cdot x_0=gN_\C \cdot x_0\subset X_\C $$
(use $M_\C \subset H_\C$).
\par   Similarly,  for $z\in X_\C$ we set
\begin{equation}\label{eq=S} S(z)=\pi_1(\pi_2^{-1}(z))\ .\end{equation}
If $z=g\cdot x_0$ for $g\in G_\C$, then notice
$S(z)=gH_\C\cdot \xi_0$.
Moreover, for $z\in X_\C$ and $\xi\in \Xi$ one has
the incidence relations
\begin{equation}\label{eq=indi}
z\in E(\xi)\iff \pi_1^{-1}(\xi)\cap \pi_2^{-1}(z)\neq \emptyset
\iff \xi\in S(z)\, .\end{equation}

\par The space of horospheres on $X_\C$ shall be denoted
by ${\rm Hor}(X_\C)$, i.e.
$${\rm Hor}(X_\C)=\{ E(\xi)\mid \xi\in \Xi\}\, .$$
Our first objective is to show that $\Xi$
parameterizes ${\rm Hor}(X_\C)$:

\begin{proposition}\label{prop=ident} The map
$$E: \Xi \to {\rm Hor}(X_\C), \ \ \xi\mapsto E(\xi)$$
is a $G_\C$-equivariant bijection.
\end{proposition}

\begin{proof} Surjectivity and $G_\C$-equivariance are clear
by definition. It remains to establish
injectivity. For that write $G_\C^{E(\xi_0)}$
for the stabilizer of $E(\xi_0)$ in $G_\C$. By $G_\C$-equivariance it
is enough to show that $G_\C^{E(\xi_0)} \subseteq M_\C N_\C $.
Assume that
$g\cdot E(\xi_0)=E(\xi_0)$. Then $gN_\C\subseteq N_\C H_\C$. In particular,
$g=nh\in N_\C H_\C$. As $G_\C^{E(\xi_0)}$ is a group, and $n\in G_\C^{E(\xi_0)}$, it
follows, that $h\in G_\C^{E(\xi_0)}$. By Lemma \ref{l-12one}
from below it follows that
$h\in M_\C$. Hence $g=h(h^{-1}nh)\in M_\C N_\C$, as $M_\C$ normalizes $N_\C$.
\end{proof}

\begin{lemma}\label{l-12one} Assume that $h\in H_\C$ is such that
$h \cdot E(\xi_0)=E(\xi_0)$. Then $h\in M_\C$.
\end{lemma}

\begin{proof} Identify the tangent space $T_{x_0}(G_\C /H_\C )$ with
$\fg_\C /\fh_\C$.
Then, as $(hN_\C h^{-1})\cdot x_0=N_\C \cdot x_0$, it follows that
$$\Ad (h)(\fn_\C\oplus \fh_\C)=\fn_\C\oplus \fh_\C\, .$$
Thus, if $U\in \fn_\C$, there exists $Z\in \fn_\C$ and  $L\in \fh_\C$ such that
$\Ad (h)U = Z+L$. Applying $({\bf 1}-\tau)$ this equality, we get
$\Ad (h)(U-\tau (U))=Z-\tau (Z)$. As $\fq_\C=({\bf 1}-\tau)(\fn_\C)\oplus \ft_\C$, and
this sum is orthogonal with respect to Killing form, it follows that $\Ad (h)\ft_\C
=\ft_\C$. In particular, $h\in N_{H_\C}(\ft_\C)$.
\par We recall the Riemannian dual $X^r=G^r/K^r$ of $X=G/H$ which corresponds
to the Lie algebras $\fg^r=\fk^r + \fs^r$ with $\fk^r=(\fh\cap \fk)
+ i (\fh\cap \fs)$ and $\fs^r=i(\fq\cap \fk) + (\fq\cap\fs)$.
Notice that $\fa$ is maximal abelian in $\fs^r$.

\par To continue with the proof, we observe that $N_{H_\C}(\ft_\C)=N_{K^r}(\fa)M_\C$.
Thus we may assume that $h\in N_{K^r}(\fa)$. Write $\sigma_r$ for the complex conjugation
in $G_\C$ with respect to the real form $G^r$. Then
taking $\sigma^r$ fixed points in $hN_\C\in N_\C H_\C$ yields
$h N^r\in N^r K^r$ with $N^r=G^r\cap N_\C$. Thus the situation is
reduced to the Riemannian case where it follows from
\cite{hel}, p.78.
\end{proof}

It is crucial to observe
that there is a $T_\C$-action on $\Xi$ which commutes with the
left $G_\C$-action:

\begin{proposition} \label{prop=comm}Let $\xi=g\cdot\xi_0\in \Xi$, $g\in G_\C$.
For $t\in T_\C$ the prescription
\begin{equation} \label{eq=h-action}\xi\cdot t=gt\cdot\xi_0
\end{equation}
defines an element of $\Xi$. In particular,
\begin{equation} \label{eq=h-action2}
T_\C\times \Xi\to \Xi, \ \ (t,\xi)\mapsto \xi\cdot t
\end{equation}
defines an action of $T_\C$ on $\Xi$, which commutes with the
natural action of $G$ on $\Xi$.
\end{proposition}

\begin{proof} As  $T_\C$ normalizes $M_\C N_\C$ it follows
that (\ref{eq=h-action}) is  defined.
Finally, (\ref{eq=h-action}) implies that (\ref{eq=h-action2}).
defines a left-action of $T_\C$.
\end{proof}

It is obvious
that the map
\begin{equation}\label{eq=haction3}
(G_\C \times T_\C)\times \Xi\to \Xi, \ \ \left((g,t), \xi\right)
\mapsto g\cdot\xi\cdot t
\end{equation}
is a holomorphic action of the complex group $G_\C\times T_\C$
on the homogeneous space $\Xi$.

\par The remainder of this section will be devoted to
the definition and basic discussion of an important
$G\times T$-invariant subset
$\Xi_+$ of $\Xi$.

\par We recall from Subsection \ref{ss=11} the
polydisc $\cT_+ =TA_+$ and define
$$\Xi_+= G\cT_+\cdot \xi_0=GA_+\cdot\xi_0\, .$$
We record that $\Xi_+$ is a $(G\times T)$-invariant
subset of $\Xi$.

\par The set $G P_\C $ is open in $G_\C$ and $G\cap P_\C = MT$. Hence
$G/MT$ can be viewed as an open, complex submanifold of the flag
manifold $F=G_\C /P_\C$. We write $F_+=G P_\C/P_\C$ for the image of $G/MT$ in $F$
and call $F_+$ the {\it flag domain}. Although obvious we emphasize
that  $F_+$ is $G$-homogeneous.
\par Notice $G/MT$ is the base space  of the holomorphic fiber bundle
$G/M\times_T\cT_+\to G/MT$ with fiber $\cT_+/\Gamma$.
There is a natural action of $G\times T$
on $G/M\times_T\cT_+$ given by

$$(G\times T)\times \left(G/M\times_T\cT_+\right)\to
G/M\times_T\cT_+, \ \ \left((g,t), [xM,a]\right)
\mapsto [gxM, at]\, . $$

The next lemma gives us basic structural information on
$\Xi_+$.
\begin{lemma}\label{lemma=iso}
The set $\Xi_+$ is open in $\Xi=G_\C/ M_\C N_\C$.
Moreover, the mapping
$$\Phi: G/M\times_T \cT_+\to \Xi, \ \  [gM,t]\mapsto gt\cdot \xi_0$$
is a $G\times T$-equivariant biholomorphism onto $\Xi_+$.
\end{lemma}

\begin{proof}  Clearly,  $\Phi$ is a defined
$G\times T$ -equivariant map with $\im \Phi=\Xi_+$.
By the definition of the complex structure of
$G/MT$ the holomorphicity of the map is clear, too. Let us
show that $\Phi$ is injective.
For that assume that $g_1t_1\cdot \xi_0 = g_2t_2\cdot \xi_0$, $g_j\in G$,
$t_j\in \cT_+$. By $G$-equivariance we may assume that
$g_2={\bf 1}$. Then $g_1\in G\cap P_\C =MT$ and
w.lo.g. we may assume that $g_1\in M$.
Consequently, as $T_\C \cap M_\C N_\C=\Gamma$, we obtain $t_1\in t_2\Gamma$,
i.e. $[M,t_1]=[M,t_2]$. Hence $\Phi$ is injective.
\par A standard computation yields that $d\Phi$ is an immersion
and a simple dimension count shows that $\dim G/MT + \dim \cT_+
=\dim \Xi$. In particular, $\Phi$ is a submersion
and $\im \Phi=\Xi_+$ is open, concluding the proof
of the lemma.
\end{proof}
\subsection{Fiberings.} To conclude this section we mention three  natural
fibrations in relation to $\Xi_+$ and $F_+$.
\par Write  $S^+= \exp (\fs^+)$  and recall  that
the map
$$Y=G/K\to G_\C /K_\C S^+, \ \ gK\mapsto gK_\C S^+$$
is a $G$-equivariant open embedding. Henceforth
$Y$ will be understood as an open subset of
the flag manifold $G_\C/ K_\C S^+$.

\begin{lemma} The following assertions hold:
\begin{enumerate}
\item The natural map
$$\Xi_+\to F_+, \ \ z M_\C N_\C \mapsto z P_\C$$
is a holomorphic fibration with fiber
${\mathcal T}_+/\Gamma$.
\item  The natural map
$$F_+\to Y, \ \ gMT\mapsto gK$$
is a holomorphic fibration with fiber the flag variety $K/MT$.
\item The natural map
$$\Xi_+\to Y, \ \ gt\cdot\xi_0\mapsto gK$$
is a holomorphic fibration  with fiber
$K/M\times_T \cT_+$.
\end{enumerate}
\end{lemma}

\begin{proof} (i) follows from
$G\cap P_\C = MT$ and (ii) is obvious. Finally (iii) is a consequence
(i) and (ii).
\end{proof}

\section{The $G\times T$-Fr\'echet module $\cO(\Xi_+)$}

The natural  action of $G\times T$ on $\Xi_+$
gives rise to a representation of $G\times T$ on the Fr\'echet
space  $\cO(\Xi_+)$ of holomorphic functions on $\Xi_+$.
We will decompose $\cO(\Xi_+)$ with respect to this
action.  By the compactness of $T$,
it is clear that $\cO(\Xi_+)$ decomposes discretely
under $T$. It turns out that each $T$-isotypical component
is the section module of a holomorphic line bundle over the flag
domain $F_+$ and that all such section modules arise in this
manner.

\par In the second part of this section we turn our attention to
$G\times T$-invariant Hilbert spaces of holomorphic functions
on $\Xi_+$. By definition these are unitary $G\times T$-modules
$\cH$ with continuous $G\times T$-equivariant embeddings
into $\cO(\Xi_+)$. There are many interesting
examples such as weighted Bergman and weighted Hardy spaces.
We will discuss the Hardy space $\Hsp$ on $\Xi_+$ with constant
weight and show that $\Hsp$ constitutes a natural model for the
the $H$-spherical holomorphic discrete series of $G$.

\subsection{The decomposition of $\cO(\Xi_+)$}
In Section \ref{section=hor-I} we exhibited a natural action of
$G\times T$ on $\Xi_+$, namely

\begin{equation}\label{h4}
(G\times T)\times \Xi_+\to \Xi_+, \ \ \left((g,t),\xi\right)
\mapsto g\cdot\xi\cdot t\ .\end{equation}
We recall that $\cO(\Xi_+)$ becomes a Fr\'echet space when endowed with
the topology of compact convergence.

\begin{remark} Finite dimensional representation theory of $G_\C$
shows that $\Xi$ (and hence $\Xi_+$)
is holomorphically separable. In particular
$\cO(\Xi_+)\neq \{0\}$. \end{remark}

Denote by
${\rm GL} (\cO(\Xi_+))$ the group of bounded invertible operators
on $\cO(\Xi_+)$.
\par The action
(\ref{h4}) induces a continuous representation
of $G\times T$ on $\cO(\Xi_+)$:

$$L\otimes R: G\times T \to {\rm GL} (\cO(\Xi_+)), \ \
\left((L\otimes R)(g,t) f\right)(\xi)=f(g^{-1}\cdot \xi\cdot t^{-1})\, ,$$
$(g,t)\in G\times T$, $f\in \cO(\Xi_+)$, and $\xi\in \Xi_+$.

We first decompose $\cO(\Xi_+)$ under the action of the compact
torus $T$. Denote by
$\widehat {T/\Gamma}$ the  character group
of $T/\Gamma$, i.e. $\widehat{T/\Gamma}={\rm Hom}_{\rm cont}(T/\Gamma, {\Bbb S}^1)$.
In the sequel we identify $\widehat{T/\Gamma}$ with the lattice

$$ \Lambda=\{\lambda \in \fa^*\mid
\forall U\in (\exp|_\ft)^{-1}(\Gamma)\,\, \lambda (U)\in 2\pi i \Z\}\, .$$
Explicitly, to $\lambda\in \Lambda$ one associates the character
$\chi_\lambda(t\Gamma)=e^{\lambda(\log t)}$.
Often we will write $t^\lambda$ for $\chi_\lambda (t\Gamma)$.

\par The assumption that $G_\C$ is simply connected allows
an uncomplicated description of the lattice $\Lambda$.

\begin{lemma} $\Lambda=\left\{ \lambda\in \fa^*\mid (\forall
\alpha\in \Delta)\ {\la \lambda, \alpha\ra \over \la \alpha,\alpha\ra}
\in \Z\right\}$.
\end{lemma}

\begin{proof}  $''\subseteq''$: Let $\lambda\in \Lambda$.
We first show that
${\la \lambda, \alpha\ra \over \la \alpha,\alpha\ra}
\in \Z$ for all $\alpha\in\Delta_k$. For that observe that
the compact symmetric space $K/ H\cap K$ embeds into $G/H$ via
the natural map
$$K/H\cap K\to G/H, \ k(H\cap K)\mapsto kH\, .$$
Thus \cite{he84}, Ch. V, Th. 4.1, yields that
${\la \lambda, \alpha\ra \over \la \alpha,\alpha\ra}
\in \Z$ for all $\alpha\in\Delta_k$. To complete the
proof of $''\subseteq''$ we still have to verify
${\la \lambda, \alpha\ra \over \la \alpha,\alpha\ra}
\in \Z$ for all $\alpha\in\Delta_n$. Fix $\alpha\in\Delta_n$.
Standard structure theory implies that there is an
embedding of symmetric Lie algebras
$(\mathfrak{su}(1,1), \mathfrak{so}(1,1))\to (\fg,\fh)$ such that
$\left[\begin{matrix}  i & 0\\ 0 & -i\end{matrix}\right]\in\mathfrak{su}
(1,1)$  is mapped to $i\check\alpha\in \ft$.
As $G_\C$ is simply connected, we thus obtain an immersive map
${\rm SU}(1,1)/ {\rm SO(1,1)}\to G/H$. In particular,
${\la \lambda,\alpha\ra \over \la \alpha, \alpha\ra}\in \Z$ must hold.

\par $''\supseteq''$: Suppose that
${\la \lambda, \alpha\ra \over \la \alpha,\alpha\ra}
\in \Z$ holds for all $\alpha$. Recall the extension
$\ft\subseteq \fc$ of $\ft$ to a compact Cartan
subalgebra of $\fg$.  In the sequel we consider $\lambda$ as
an element of $\fc^*$ which is trivial on $\fc\cap \fh$.
On p. 537 in \cite{he84}, it is
shown that $\lambda$ is analytically integral for $C=\exp \fc$ (again this needs that
$G_\C$ is simply connected). In particular $\lambda$ defines an element
$\chi_\lambda\in \hat T$. It remains to show that
$\chi_\lambda|_\Gamma={\bf 1}$. As $M=Z_{H\cap K}(\fa)$ and $\Gamma=M\cap T$,
this reduces to an assertion on the compact symmetric space
$K/H\cap K$, where it follows from
\cite{he84}, Ch. V, Th. 4.1.
\end{proof}

For each $\lambda\in \Lambda$ define the $\lambda$-isotypical  component of $\cO(\Xi_+)$ by

\begin{equation}
\cO(\Xi_+)_\lambda =
\{f\in \cO(\Xi_+)\mid ( \forall t\in  T)\  R(t)f=t^{\lambda} f\}\ .
\end{equation}
As $(R, \cO(\Xi_+))$ is a continuous representation of the
compact torus $T$ on a Fr\'echet space, the Peter-Weyl theorem
yields
\begin{equation}
\cO(\Xi_+) =\bigoplus_{\lambda \in \Lambda}\cO(\Xi_+)_\lambda
\end{equation}

Each $\cO(\Xi_+)_\lambda$ is a $G$-module for the representation $L$.
In order to describe them explicitly we recall
some facts on holomorphic line bundles.

For $\lambda\in \Lambda$ we write $\C_\lambda$ for $\C$ when considered
as a $MT$-module with trivial $M$-action and $T$ acting by $\chi_\lambda$.
Recall that $G/MT$ inherits a complex manifold structure through
its identification with the flag domain $F_+$.
In particular, to each $\lambda\in \Lambda$
one associates the holomorphic line bundle
\begin{equation}\label{eq=lineb}
\cL_{\lambda}=G\times_{MT}\C_{-\lambda}\, .
\end{equation}
Write $\cO(\cL_{\lambda})$ for its $G$-module of holomorphic sections, i.e.
 $\cO(\cL_{\lambda})$ consists of smooth functions
$f: G\to \C$ such that
\begin{itemize}\item $f(gmt)=t^{-\lambda} f(g)$ for $g\in G$, $t\in T$
and $m\in m$.
\item $G/MT\to\cL_{\lambda}, \ \ gMT\mapsto [gMT, f(g)]$ is
holomorphic.
\end{itemize}

\par The restriction of $\cL_\lambda$ to
the flag variety $K/MT$ yields the holomorphic line
bundle
$$\cK_\lambda=K\times_{MT} \C_{-\lambda}$$
over $K/MT$.
Write $\Lambda_0$ for the
$\Delta_k^-$-dominant elements of
$\Lambda$, i.e.

\begin{equation}
\Lambda_0=\left\{\lambda\in \Lambda\mid (\forall \alpha\in \Delta_k^+)
\ \la \lambda, \alpha\ra \leq 0\right\}\, .
\end{equation}

According to Bott \cite{b}, $V_\lambda=\cO(\cK_\lambda)$
is of finite dimension,  and non-trivial if and only if $\lambda\in\Lambda_0$.
By
$\cL_\lambda=G\times_{MT}\C_{-\lambda}\simeq
G\times_K (K\times_{MT} \C_{-\lambda})
$
we retrieve the standard isomorphism
$$\cO(\cL_\lambda)\simeq \cO(G\times_K V_\lambda)\, .$$
In particular,
\begin{equation}\label{eq=standard}
\cO(\cL_\lambda)\neq\{0\} \quad\iff\quad \lambda\in \Lambda_0\, .
\end{equation}

\par We remind the reader that the
$T$-weight spectrum of $\pi_\lambda$ is
contained in $\lambda+\Z_{\geq 0} [\Delta^+]$.
In particular, $\cO(\cL_{\lambda})$, if irreducible,  is a lowest
weight module for $G$
with respect to the positive system $\Delta^+$ and lowest  weight $\lambda$.

\par Finally we establish the connection between $\cO(\Xi_+)_\lambda$
and $\cO(\cL_\lambda)$. For that let us denote by $\Xi_0$
the pre-image of $F_+$ in $\Xi$, i.e.
$$\Xi_0=GT_\C\cdot\xi_0\, .$$
Notice that $\Xi_+\subset \Xi_0$.
Holomorphicity and $T$-equivariance yield $\cO(\Xi_+)_\lambda
=\cO(\Xi_0)_\lambda$. Likewise holds for $\cO(\cL
_\lambda)$.
Thus holomorphic extension and restriction gives a natural
$G$-isomorphism $\cO(\Xi_+)_\lambda\simeq \cO(\cL_\lambda)$.
\par We summarize our discussion.

\begin{proposition} The $G\times T$-Fr\'echet module $\cO(\Xi_+)$
decomposes as
$$\cO(\Xi_+)=\bigoplus_{\lambda\in \Lambda_0} \cO(\Xi_+)_\lambda\ .$$
Moreover, holomorphic extension and restriction canonically
identifies $\cO(\Xi_+)_\lambda$
with the section module $\cO(\cL_\lambda)$.
\end{proposition}
We conclude this subsection with some comments
on unitarization of the section modules $\cO(\cL_\lambda)$.

\begin{remark} Let $\lambda\in \Lambda_0$ and let us denote by $\cO(\cL_\lambda)^{K-{\rm fin}}$  the
$(\fg,K)$-module of $K$-finite sections of $\cO(\cL_\lambda)$.
Let us assume that $\cO(\cL_\lambda)^{K-{\rm fin}}$ is
irreducible. Then $\cO(\cL_\lambda)^{K-{\rm fin}}$
identifies with the generalized Verma module
$N(\lambda)=\cU(\fg_\C)\otimes_{\cU(\fk_\C +\fs^-)} V_\lambda$
and the Shapovalov form on $N(\lambda)$ gives rise
to the (up to scalar unique) contravariant Hermitian form on
$\cO(\cL_\lambda)^{K-{\rm fin}}$. We say that
$\cO(\cL_\lambda)^{K-{\rm fin}}$ is {\it unitarizable} if the
Shapovalov form is positive definite. Another way to formulate  it
is that there exists a unitary lowest weight representation
$(\pi_\lambda, {\mathcal H}_\lambda)$ such that the $(\fg, K)$-module
of $K$-finite vectors
${\mathcal H}_\lambda^{K-{\rm fin}}$ is $(\fg,K)$-isomorphic to
$\cO(\cL_\lambda)^{K-{\rm fin}}$.  In this situation
$\cO(\cL_\lambda)$ is then naturally $G$-isomorphic to
the hyperfunction vectors ${\mathcal H}_\lambda^{-\omega}$
of $\pi_\lambda$.
\par We want to emphasize that not all $\lambda\in \Lambda_0$
correspond to unitarizable modules $\cO(\cL_\lambda)^{K-{\rm fin}}$
(a necessary condition is $\lambda|_\Omega\geq 0$ and we
refer to \cite{EHW} for more precise information).
However, we want to stress that
$\cO(\cL_\lambda)^{K-{\rm fin}}$ is automatically unitarizable
if $\lambda|_\Omega$ is sufficiently positive
(for example if condition (\ref{eq=HC})  below is
satisfied).
\end{remark}

\subsection{The Hardy space on $\Xi_+$}

The objective of this section is to introduce
the Hardy space on $\Xi_+$ and to  prove some of its basic
properties.

\par We begin with some measure theoretic preliminaries.
The groups $G_\C$ and $M_\C N_\C$ are unimodular, and hence $\Xi=G_\C /M_\C N_\C$
carries a $G_\C$-invariant measure $\mu$.
\par Recall that $M$ is a compact subgroup of $G$ and denote by
$dm$ a normalized Haar measure on $M$. Further we let $dg$ and $d(gM)$
denote left $G$-invariant measures on $G$, resp. $G/M$, normalized
subject to the condition
$$\int_G f(g)\ dg=\int_{G/M}\int_M f(gm)\ dm\ d(gM)$$
for all $f\in L^1(G)$.
\par Notice that
the stabilizer in $G$ of any point $\xi\in \cT_+\cdot\xi_0\subset \Xi_+$ is the
compact subgroup $M$. In particular one has

\begin{equation}\label{eq-integral}
\int_G f(g\cdot \xi)\, dg
=\int_{G/M} f(g\cdot \xi)\, d(gM)
\end{equation}
for all $\xi\in \cT_+\cdot\xi_0$ and integrable functions $f$ on $\Xi_+$.

\par Write $\|\cdot \|_2$ for the
$L^2$-norm on $L^2(G)$. Let us
remark that the representation $(R, \cO(\Xi_+))$ of $T$ naturally
extends to a representation of the semigroup
$t\in  \cT_-\cup T $, also denoted by $R$.
Furthermore if $f\in \cO(\Xi_+)$ and $t\in \cT_-$ then we can
define the restriction of $R(t)f$ to $G$ by
$R(t)f|_G : G\to \C$ by $R(t)f|_G (g) =
f(gt^{-1}\cdot \xi_0)$. The Hardy norm of $f\in \cO(\Xi_+)$ is
defined by
\begin{equation}\label{def-Hardynorm}
\|f\|^2=\sup_{t\in  \cT_+}\int_G | f(gt\cdot\xi_0)|^2\, dg
=\sup_{t\in  \cT_-}\| R(t)f|_G\|^2_2\, .
\end{equation}
Let
\begin{equation}\label{de-Hardy}
\cH^2 (\Xi_+ )= \{f\in \cO (\Xi_+)\mid \| f\|<\infty \}\, .
\end{equation}
Obviously
\begin{equation}\label{eq=contract}
\|R(t)f\| \le \|f\| \qquad \hbox{for all $t\in \cT_-$} \end{equation}
and hence $\cT_-$ acts on $\Hsp$ by contractions.
Note, that $R(t)f|_G$ is right $M$-invariant, and, by the definition
of the Hardy space norm $R(t)f|_G\in L^2(G/M)\subseteq L^2(G)$.

\begin{lemma} The space $\Hsp$ is a Hilbert space. Furthermore, the
following holds:
\begin{enumerate}
\item
For $\xi\in \Xi_+$ the point evaluation map
 $\ev_\xi :\Hsp\ni f\mapsto f(\xi )\in \C$
is continuous.
\item The boundary value map $\beta :\Hsp \to L^2(G/M)\subseteq L^2(G)$
$$\beta (f)=\lim_{\cT_-\ni t\to e}R(t)f|_G$$
is an isometry into $L^2(G/M)$.
\end{enumerate}
\end{lemma}
\begin{proof} The proof follows a standard procedure
and will be more a sketch. We refer to
\cite{hoo91}, in particular the proof of Theorem 2.2, for a
detailed  discussion of the underlying methods.

\par Let $\xi \in \Xi_+$. Then there exist relatively compact
open sets $U_G\subseteq G$ and $U_T\subseteq \cT_+$ such that
$\xi \in U_GU_T\cdot \xi_0$. Thus, there is a constant $c>0$ such that
the Bergman-type estimate
$$\int_{U_GU_T\cdot\xi_0}|f(\xi )|^2\, d\mu(\xi)\le c\cdot \| f\|^2\, $$
holds for all $f\in \Hsp$.
This implies that $\Hsp $ is complete, and that
point evaluations are continuous.
\par Write $\C_+=\{z\in \C\mid \mathrm{Im}(z)>0\}$
for the upper half plane and fix $Z\in i\gO$. We
notice that the map $\cT_-\ni t\mapsto R(t)f|_G\in L^2(G)$
is well defined and holomorphic.
Hence
$$L_f : \C_+\to L^2(G/M); \ L_f(z)= R(\exp (z Z))f|_G\in L^2(G/M)$$
defines a holomorphic function on $\C_+$.

By Lemma 2.3 in \cite{hoo91} it follows, that
$\lim_{z\to 0}L_f(z)$ exists, and is monotonically
increasing as $s\searrow 0$ along each line segment
$\exp (siZ)$, or, because of the right invariance
of $dg$, on each $t\exp (siZ)$, $t\in  T$. As in \cite{hoo91}, one shows, that this limit
is independent of $Z$. Thus, we get a boundary value map
$\beta : \Hsp \to L^2(G/M)$, defined by
$$\beta (f)=\lim_{t\to e}R(t)f|_G\, .$$
By the definition of the Hardy space norm, we obviously have
$$\|\beta (f)\|_2\le \|f\|\, .$$
But, as the norm $\| R(\exp (sZ))f\|_2$ is monotonically increasing for
$s\searrow 0$, it follows that
$$\|R(\exp (sZ))f\|_2\le \|\beta (f)\|$$
for all $s\in \R^+$. Thus
$$\|R(t)f|_G\|\le \|\beta (f)\|\, .$$
It follows, that $\beta : \Hsp \to L^2(G)$ is an isometry, and
hence $\Hsp $ is a Hilbert space.
\end{proof}

Clearly $L\otimes R$ defines a unitary representation of
$G\times T$ on $\Hsp$. We are going to decompose
$\Hsp$ with respect to this action. As before
we begin with the decomposition under $T$.
For $\lambda\in \Lambda$ the
$\lambda$-isotypical  component of $\Hsp$
is given by $\Hsp_\lambda=\Hsp\cap \cO(\Xi_+)_\lambda$.
The Peter-Weyl theorem yields the orthogonal
decomposition

\begin{equation}
\Hsp =\bigoplus_{\lambda \in \Lambda_0}\Hsp_\lambda\,
\end{equation}
of $\Hsp$ in $G$-modules.

We draw our attention to the unitary $G$-modules
$\Hsp_\lambda$ inside of $\cO(\Xi_+)_\lambda$.

\par Suppose that $\Hsp_\lambda\neq \{0\}$. Then
$\cO(\cL_{\lambda})\neq \{0\}$ and the restriction mapping
$$\Hsp_\lambda\to \cO(\cL_{\lambda})$$
gives a $G$-equivariant embedding.
Moreover $\beta(\Hsp_\lambda)\subset L^2(G)$.
Thus $\Hsp_\lambda$ is a module of the holomorphic discrete series of $G$.
In terms of $\lambda$ this means that
$\lambda$ satisfies the Harish-Chandra condition \cite{hc55}
\begin{equation}\label{eq=HC}
\langle \lambda -\rho(\fc),\alpha\rangle >0 \qquad (\forall \alpha\in \Sigma_n^+)\, ,
\end{equation}
where $\rho(\fc)={1\over 2}
\sum_{\alpha\in \Sigma^+}\alpha$.

Write $\Lambda_{{\rm sd}}$ for the set of all $\lambda\in \Lambda_0$ which satisfy
(\ref{eq=HC}).

\par Conversely, let $\lambda\in \Lambda_{\rm sd}$ and write ${\mathcal H}_\lambda$ for a
corresponding unitary lowest weight module with lowest  weight $\lambda$.
Denote by $v_\lambda\in {\mathcal H}_\lambda$ a normalized lowest   weight vector and
write $d(\lambda)$ for the formal dimension
(see \cite{hc55} or (\ref{***}) below).
It is then straightforward
that
$${\mathcal H}_\lambda\to \Hsp, \ \ v\mapsto \left(gt\cdot\xi_0\mapsto \sqrt{d(\lambda)}\cdot
t^{-\lambda}\langle \pi_\lambda(g^{-1})v, v_\lambda\rangle\right) $$
defines a $G$-equivariant isometric embedding.  Hence $\Hsp_\lambda
\simeq {\mathcal H}_\lambda\neq \{0\}$.

Summarizing our discussion we obtain the Plancherel decomposition  for $\Hsp$.

\begin{proposition}\label{th=Plancherel} As a $G$-module
the Hardy space decomposes
as
$$\Hsp \simeq \bigoplus_{\lambda\in \Lambda_{\rm sd}} {\mathcal H}_\lambda\, .$$
\end{proposition}

\begin{remark}\label{rem=disc} (a) The set $\Lambda_{\rm sd}$ describes
the set of all
$H$-spherical unitary lowest weight representations (up to equivalence)
whose matrix coefficients are square integrable on $G$, i.e.
$\Lambda_{\rm sd}$ is the spectrum of the $H$-spherical holomorphic
discrete series of $G$.
\par  \noindent(b) Later we well mainly deal
with the spectrum $\Lambda_2$ of the holomorphic
discrete series on $X$. One has
$$\Lambda_2\subseteq \Lambda_{\rm sd}$$
with equality precisely for the equal rank cases \cite{oo88,oo91}.
\end{remark}

\section{Complex Horospheres II: Horospheres with no real points}
\label{section=hor-II}
\noindent
We continue our discussion of complex horospheres from Section \ref{section=hor-I}.
We will introduce the notion of horosphere without real points and investigate
$\Xi_+$ with respect to this property. In addition we will
prove some dual statements for the minimal tubes $D_{\pm}$.

\begin{definition} We say that the complex horosphere $E(\xi)\subset X_\C$
 has
\textit{no real points} if $E(\xi)\cap X =\emptyset$. We denote by
$\Xi_{nr}\subset \Xi$ the subset of those $\xi$ which correspond to
horospheres with no real points.
\end{definition}

\begin{lemma} The set $\Xi_{nr}$ is a $G$-invariant
subset of $\Xi$.
\end{lemma}

\begin{proof} Let $\xi \in \Xi_{nr}$ and $g\in G$. Assume, that
$ x\in E(g\cdot \xi)\cap X$. Then $g^{-1}x\in E(\xi) \cap X$, contradicting the
assumption that $E(\xi)$ has no real points.
\end{proof}

Recall the open $G$-invariant subset $\Xi_+=GA_+\cdot\xi_0\subset \Xi$.
In the sequel it will be useful to consider with $\Xi_+$ its
pre-image $\widetilde \Xi_+$ in $G_\C$, i.e.
$$\widetilde \Xi_+=GA_+ M_\C N_\C\ .$$
It is clear that $\widetilde\Xi_+$ is a left $G$ and right $M_\C N_\C$ invariant
open subset of $G_\C$.

\par Next we draw our attention to the Zariski open subset
$N_\C A_\C H_\C  $ of $G_\C$. Our  objective
is to study $\widetilde \Xi_+$ in relation to $N_\C A_\C H_\C $.

\begin{remark} Notice that $\widetilde\Xi_+^{-1}\subset N_\C A_\C H_\C $
is equivalent to $G\subset N_\C A_\C H_\C $. However the latter is
true only for ${\rm rank}\, X=1$, i.e. $\dim \ft =1$.  In general,
$G\cap N_\C A_\C H_\C $ is an open and dense subset of $G$
(cf. Theorem \ref{th=monoton} below).
\end{remark}

There is a right $H_\C$ and left $N_\C$-invariant
holomorphic middle-projection
$$a_H: N_\C A_\C H_\C  \to A_\C/ \Gamma, \ \ z\mapsto a_H(x)$$
In particular, for each $\lambda\in \Lambda$ we obtain natural
$(N_\C, H_\C)$-invariant holomorphic maps
$$N_\C A_\C H_\C  \to \C, \ \ x\mapsto a_H(x)^\lambda\ .$$
The holomorphic function $N_\C A_\C\cdot x_0
\to A_\C/\Gamma$ induced by $a_H$ shall also be denoted by
$a_H$.

\par The function $a_H$ enables us to give a
useful geometric description of horospheres.

\begin{lemma}\label{lemma=Hc} Let $\xi=g\cdot\xi_0\in \Xi$ for
$g\in G_\C$. Then
\begin{eqnarray*} E(\xi) &=&
\{ z\in X_\C \mid g^{-1}z\in N_\C A_\C\cdot x_0, \ a_H(g^{-1}z)=\Gamma\}\\
&=&\{ z\in X_\C \mid g^{-1}z\in N_\C A_\C\cdot x_0, \ a_H(g^{-1}z)^\lambda=1
\quad\text{for all $\lambda\in \Lambda$}\}
\end{eqnarray*}
\end{lemma}

\begin{proof} $''\subseteq''$: If $z\in E(\xi)$, then $z=gn\cdot x_o$
for some $n\in N_\C$. Thus $g^{-1}z=n\cdot \xi_0\in N_\C A_\C \cdot x_0$
and $a_H(g^{-1}z)=a_H(n\cdot x_0)=\Gamma$.
\par $''\supseteq''$: Conversely, let $z\in X_\C$ be such that
$g^{-1}z\in N_\C A_\C \cdot x_0$ and $a_H(g^{-1}z)=\Gamma$.
>From the first condition follows that $g^{-1} z =na\cdot x_0$
for some $n\in N_\C$ and $a\in A_\C$; the second condition implies
$a\in \Gamma$. Thus $z\in g\cdot\xi_0$, as was to be shown.
\end{proof}

We define a subset of $\Lambda_0$ by
\begin{eqnarray}\label{lattice}
\Lambda_{\geq 0}&=&\{\lambda\in \Lambda_0\mid \lambda|_\Omega\geq 0\}\\
&=&\left\{\lambda\in \Lambda\mid \lambda|_\Omega\geq 0, \
(\forall \alpha\in \Delta_k^+)\quad \la \lambda, \alpha\ra \leq 0\right\}\, .
\end{eqnarray}
The following theorem is the main geometric result of the paper.

\begin{theorem}\label{th=monoton}The following assertions hold:
\begin{enumerate}
\item $G\cap N_\C A_\C H_\C $ is  open and dense
in $G$.
\item Let $\lambda\in\Lambda_{\geq 0}$. Then the function $a_H^{\lambda}|_{G\cap N_\C A_\C H_\C }$
extends to a continuous function on $G$ and
$$|a_H(g)^\lambda|\leq 1\qquad(g\in G) .$$
\end{enumerate}
\end{theorem}

\begin{proof} The approach to prove this theorem lies in the use of the structural
decomposition
\begin{equation}\label{decomp} G=KA_qH \end{equation}
where $A_q=\exp(\fa_q)$ with $\fa_q\subseteq \fs\cap \fq$
a maximal abelian subspace. There is a natural way to construct
a flat $\fa_q$ out of the weight space decomposition $\fg_\C=\fa_\C +\fm_\C + \bigoplus_{\alpha\in \Delta}\fg_\C^\alpha$.
It will be briefly reviewed.
Let $\gamma_1,\ldots ,\gamma_r\in \gD^+_n$ be a maximal set of long strongly orthogonal roots.
Then  one can find
$Z_j\in \fg_\C^{\gamma_j}$, $j=1,\ldots ,r$, such
that
\begin{equation}\label{def-ap}
\fa_q=\bigoplus_{j=1}^r\R (Z_j-\tau (Z_j))
\end{equation}
is a maximal abelian subspace of $\fs\cap \fq$;
further
\begin{equation} \label{inc} A_q\subset S^+  \oline{A_-} H_\C \end{equation}
(see \cite{ho}, p. 210-211 for all that).
\par(i)  As $S^+\subseteq N_\C$ , we obtain from (\ref{decomp}) and (\ref{inc}) that
$$G\subset N_\C K \oline{A_-}  H_\C\ . $$
Hence it is sufficient to show that

\begin{equation} \label{dense} K \oline{A_-}  \cap  N_\C^+  A_\C (H_\C \cap K_\C)
\quad\text{
is open and dense in $K \oline{A_-}$} \end{equation}

\par To continue, we first have to recall some facts related to
the Iwasawa decomposition of $K_\C$.
Write $\widetilde N_\C^+$ for a maximal $C$-stable
unipotent subgroup of $K_\C$ containing $N_\C^+$ and set $\widetilde A=\exp (i\fc)$.
Then $K_\C=\widetilde N_\C^+  \widetilde A K$ is an Iwasawa decomposition
of $K_\C$. We recall that $\Omega$ and hence
$\oline{A_-}$ is $\cW_k$-invariant. Thus
Kostant's non-linear convexity theorem (cf. \cite{he84}, Ch. IV, Th. 10.5)
implies that $K \oline{A_-} \subset \tilde N_\C^+ \oline{A_-} K$.
As $\widetilde N_\C^+\subset N_\C^+ M_\C$ and
$A\subseteq \widetilde A\subseteq AM_\C$, we thus get
$K \oline{A_-} \subset N_\C^+ \oline{A_-} M_\C K$.
In particular, in order to establish (\ref{dense}) it is enough to
verify that $K\cap N_\C^+ A_\C (H_\C \cap K_\C) $ is dense
in $K$. But this is known (for example it follows from
Lemme 2.1 in \cite{Cl88}).
\par(ii) In the proof of (i) we have seen that $G\subset N_\C M_\C \oline{A_-} K H_\C  $.
Thus we only have to show that $a_H^\lambda$ can be defined as a holomorphic
function on $K_\C$ with $|a_H(ak)^\lambda|\leq 1 $ for all $k\in K$ and $a\in \oline{A_-}$.
For that let $(\tau_\lambda, V_\lambda)$
denote the holomorphic $(H_\C\cap K_\C)$-spherical representation
of $K_\C$ with lowest weight $\lambda$. Write $(\cdot, \cdot)$ for
a $K$-invariant inner product on $V_\lambda$. Let
$v_\lambda$ be a normalized
lowest weight vector and $v_H$ be the spherical vector with
$(v_H, v_\lambda)=1$.
Then for all $x\in N_\C^+ A_\C  (H_\C\cap K_\C) \subset K_\C$ we have
$$(\pi_\lambda(x) v_H, v_\lambda)=a_H(x)^{\lambda}\ .$$
 As the left hand side has a holomorphic extension to $K_\C$, the same holds
for $a_H^{\lambda}$. Finally, for $a\in \oline{A_-}$ and $k\in K$ we have
$$a_H(ak)^{\lambda}=a^{\lambda} a_H(k)^{\lambda}\ .$$
Observe that $a^{\lambda}\leq 1 $ as $\lambda\in \Lambda_{\geq 0}$
and that $|a_H(k)^{\lambda}|\leq 1$ for all $k\in K$ by Lemma 2.3
in
\cite{Cl88}.
This completes the proof of (ii).
\end{proof}

Theorem \ref{th=monoton} features interesting and important corollaries.

\begin{corollary}\label{cor=1} Let $\lambda\in \Lambda_{\geq 0}$ be such that
$\lambda|_\Omega>0$.  Then $a_H^{\lambda}|_{\widetilde\Xi_+^{-1}\cap N_\C  A_\C H_\C }$
extends to a holomorphic function on $\widetilde\Xi_+^{-1}$ with
$$|a_H(x)^{\lambda}|< 1 \qquad (x\in \widetilde \Xi_+^{-1})\, .$$
\end{corollary}

\begin{corollary}\label{cor=Hc} $\Xi_+\subseteq
\Xi_{nr}$, i.e. $E(\xi)\cap X=\emptyset$ for all $\xi\in \Xi_+$.
\end{corollary}
\begin{proof} Suppose that there exists $\xi\in \Xi_+$ such that $E(\xi)\cap X\neq\emptyset$.
In other words $\widetilde\Xi_+\cap H_\C \neq \emptyset$
$\iff$ $\widetilde\Xi_+^{-1}\cap H_\C \neq \emptyset$
; a contradiction to
the previous corollary.
\end{proof}

\begin{remark} (Monotonicity/Convexity) Theorem \ref{th=monoton} (ii) has a natural
interpretation in terms of convexity/monotonicity. Write $\pr_{\fa}=\Im \log a_H$
and note that $\pr_\fa: N_\C A_\C H_\C \to \fa$ is a well
defined continuous map. Theorem \ref{th=monoton} (ii) is then equivalent to
the inclusion
\begin{equation} \label{eq=incl}
\pr_\fa(G\cap N_\C A_\C H_\C )\subseteq\bigoplus_{\alpha\in \Delta_n^-\cup \Delta_k^+}
\R_{\geq 0}\cdot  \check\alpha\, .
\end{equation}
\end{remark}

\subsection{Dual statements for the minimal tubes}

Recall from Subsection \ref{ss=11} the minimal tubes
$D_\pm=GA_\pm\cdot x_0$ in $X_\C$ with edge $X$.

\par It follows from Neeb's non-linear convexity theorem \cite{n94}
that
\begin{equation}\label{conv} GA_- \subseteq N_\C M_\C A_-  G\, .\end{equation}
This fact combined with Theorem \ref{th=monoton} yields
\begin{equation}\label{eq=ci}
GA_-H_\C \cap N_\C A_\C H_\C \subseteq N_\C TA_-
\exp\left(\bigoplus_{\alpha\in \Delta_k^+}
\R_{\geq 0} \cdot  \check\alpha\right) H_\C \ .\end{equation}
We have shown:

\begin{corollary}\label{cor=D} Let $\lambda\in\Lambda_{\geq 0}$ be such that
$\lambda|_\Omega>0$. Then, $a_H^\lambda|_{D_-\cap N_\C A_\C \cdot x_0}$
extends to a holomorphic function on $D_-$ such that
$$|a_H(x)^\lambda|<1 \qquad (x\in D_-)\, .$$
\end{corollary}

We recall the definition of the orbits  $S(z)\subset \Xi$
for $z\in X_\C$  (cf. equation \ref{eq=S}). The convexity inclusion
(\ref{eq=ci}) delivers the dual statement
to Corollary \ref{cor=Hc}:

\begin{corollary}\label{cor=S} $S(z)\cap G/M=\emptyset$ for all
$z\in D_-$.
\end{corollary}
\begin{proof} Let $z=ga\cdot x_0$ for $g\in G$ and $a\in A_-$.
Suppose that $S(z)\cap G/M\neq\emptyset$. As $S(z)=gaH_\C \cdot \xi_0$, this
is equivalent to $aH_\C N_\C \cap G\neq \emptyset$. In other words
$Ga\cap N_\C H_\C \neq \emptyset$; a contradiction to (\ref{eq=ci}).
\end{proof}

\begin{remark} Note that (\ref{conv}) is equivalent
to $A_+G  \subseteq G A_+M_\C N_\C$. This inclusion
exhibits interesting additional structure of $\Xi_+$;
it implies

\begin{equation}\label{newxi}
\Xi_+= GA_+G\cdot \xi_0\ .\end{equation}
\end{remark}

\begin{remark}\label{rem=oc}(Generalization to other cones)
Let $\tilde \Omega$
be a $\cW_k$-invariant convex open sharp cone in $\fa$
containing $\Omega$. A particular interesting example
is the maximal cone (denoted by $c_{\rm max}$ in
\cite{ho}). In this context we would like to mention
that the results in this section remain true for $\Omega$
replaced by $\tilde \Omega$, the obvious adjustment
of $\Lambda_{\geq 0}$ understood.
\end{remark}

\section{The horospherical Cauchy transform}
\noindent
Our geometric results from Section 4 enable us to define a
natural horospherical Cauchy kernel on $\Xi_+$. The kernel
gives rise to the horospherical Cauchy transform
$L^1(X)\to \cO(\Xi_+)$. The main result is a geometric inversion
formula for the horospherical Cauchy transform for functions in the
holomorphic discrete series on $X$.

\subsection{The horospherical Cauchy kernel}
In this subsection we define the horospherical Cauchy kernel
and the corresponding horospherical Cauchy transform.
We will introduce the
holomorphic spherical Fourier transform and relate it
the horospherical Cauchy transform.

\par To begin with
we have to  recall some features of the root system $\Delta$.
Let us denote by
$$\Pi=\{ \alpha_1, \ldots, \alpha_m\}$$
a basis of $\Delta$ corresponding to the positive system $\Delta_n^+\cup \Delta_k^-$. As ${\rm Spec}\,  \ad(Z_0)=\{ -1, 0,1\}$, it  follows
that exactly one member of $\Pi$ is non-compact, say
$\alpha_m$. Define weights $\omega_1, \ldots, \omega_m\in \fa^*$
by

$${\langle \omega_i,\alpha_j\rangle\over \langle \alpha_j, \alpha_j\rangle} =
\delta_{ij}\qquad (1\leq i, j\leq n)\, .$$
Set
$$\Lambda_{>0}=\Z_{\geq 0}\cdot \omega_1 +\ldots+ \Z_{\geq 0}\cdot
\omega_{m-1} +\Z_{>0} \cdot \omega_m\ .$$
Recall the definition of $\Lambda_{\geq 0}$ from (\ref{lattice}).

\begin{lemma} \label{lemma=pos}The following assertions hold:
\begin{enumerate}
\item $\omega_i|_\Omega>0$ for all $1\leq i\leq n$. In particular,
$\lambda|_\Omega>0$ for all $\lambda\in \Lambda_{>0}$.
\item $\Lambda_{\geq 0}=\Z_{\geq 0}\cdot \omega_1 +\ldots +\Z_{\geq 0}
\cdot \omega_m$. In particular, $\Lambda_{>0}\subset
\Lambda_{\geq 0}$.
\end{enumerate}
\end{lemma}

\begin{proof} (i) Fix $x\in \Omega$. Then $x=\sum_{\alpha\in \Delta_n^+} k_\alpha \check\alpha$
with $k_\alpha>0$. Now each $\alpha\in \Delta_n^+$ can be uniquely expressed as
$\alpha=\alpha_m +\gamma$ with $\gamma\in \Z_{\geq 0}[\Delta_k^-]$. Moreover if $\alpha=\beta$
is the highest root, then $\gamma\in \Z_{>0}[\alpha_1, \ldots, \alpha_{n-1}]$.
As $k_\beta>0$, the assertion follows.

\par (ii) Set $\Lambda_{\geq 0}'=\Z_{\geq 0}\cdot \omega_1 +\ldots +\Z_{\geq 0}
\cdot \omega_m$. We first show that $\Lambda_{\geq 0}'\subseteq \Lambda_{\geq 0}$.
For that let $\lambda\in \Lambda_{\geq 0}'$, say $\lambda=\sum_{i=1}^m k_i \omega_i$
with $k_i\in \Z_{\geq 0}$. As $\alpha_1, \ldots, \alpha_{n-1}$ constitutes
a basis of $\Delta_k^-$, it follows that ${\langle \lambda, \alpha\rangle\over
\langle \alpha, \alpha\rangle}\in \Z_{\leq 0}$ for all $\alpha\in \Delta_k^+$.
Furthermore, $\lambda|_{\Omega}\geq 0$ by (i). Hence $\Lambda_{\geq 0}'\subseteq \Lambda_{\geq 0}$.
\par Finally we establish $\Lambda_{\geq 0}\subseteq \Lambda_{\geq 0}'$. For that
fix $\lambda\in \Lambda_{\geq 0}$. Then $\lambda=\sum_{i=1}^m k_i \omega_i$ with
some real numbers $k_i$. We have to show that $k_i\in \Z_{\geq 0}$.
Now $\lambda\in \Lambda_{\geq 0}$ means in particular that
${\langle \lambda, \alpha\rangle\over
\langle \alpha, \alpha\rangle}\in \Z_{\leq 0}$ for all $\alpha\in \Delta_k^+$.
Hence ${\langle \lambda, \alpha_i\rangle\over
\langle \alpha_i, \alpha_i\rangle}\in \Z_{\geq 0}$ for all $1\leq i\leq n-1$.
It remains to show that  ${\langle \lambda, \alpha_m\rangle\over
\langle \alpha_m, \alpha_m\rangle}\in \Z_{\geq 0}$. Integrality is
clear. Also since $\R_{\geq 0}\cdot \check{\alpha_m}$
constitutes a boundary ray of the cone $\Omega$, non-negativity
follows.
\end{proof}

Define the {\it horospherical
Cauchy kernel} on $\Xi_+$ as the function
$$\cK(\xi)={1\over a_H(\xi^{-1})^{-\omega_m} -1}\cdot
\prod_{j=1}^{m-1} {1\over
1 -a_H(\xi^{-1})^{\omega_j}} \qquad
(\xi\in \Xi_+)\ .$$
In view of Corollary \ref{cor=1} and Lemma \ref{lemma=pos}(i),
the function $\cK$ is  holomorphic, left $H$-invariant
and bounded on subsets of the form $GU\cdot\xi_0$ for $U\subset A_+$
compact.

 This allows us to define for a function $f\in L^1(X)$ its
{\it horospherical Cauchy transform} by

$$\widehat{f} (\xi) =\int_{X} f(x) \cdot \cK(x^{-1}\xi)\, dx
\qquad (\xi\in \Xi_+)\, .$$
We notice that the horospherical Cauchy transform is a
$G$-equivariant continuous map
$$L^1(X)\to \cO(\Xi_+), \ \ f\mapsto \widehat{f}\, .$$

\begin{remark} (a) The horospherical Cauchy kernel $\cK$ is
tied to the geometry of the minimal cone $\Omega$:
there is no larger $\cW_K$-invariant open convex cone
$\tilde\Omega$ such that $\cK$ would be holomorphic
on $G\exp(\tilde\Omega)\cdot \xi_0$ (this follows from
Lemma \ref{lemma=pos} and (\ref{eq=br})). In this context we
wish to point
the difference to the results of Section \ref{section=hor-II}
which are valid for a wider class of convex cones (cf. Remark
\ref{rem=oc}).
\par  \noindent(b) For each $\lambda\in \Lambda_0$ and $\xi\in \Xi_+$
consider the complex hypersurface
$$L(\lambda, \xi)=\{z\in X_\C\mid {a_H(\xi^{-1} z)^\lambda -1=0\}}$$
in $X_\C$. Their intersection is $E(\xi)$
and they do not  intersect $X$.
The singular set of the
horospherical Cauchy kernel is the union of the  $m$ hypersurfaces
$L(\omega_i, \lambda)$ and the
edge of this set is just $E(\xi)$. It means that if $f$
is boundary value of a holomorphic
function on $D_+$ then $\widehat f$ is a residue on $E(\xi)$.
\end{remark}

\par The horospherical Cauchy transform can be decomposed in its constituents
associated to the elements $\lambda\in \Lambda_{>0}$. More precisely,
for $\lambda\in \Lambda_{>0}$ and $f\in L^1(X)$ let us define
$\widehat{f}_\lambda\in \cO(\Xi_+)$ by
$$\widehat{f}_\lambda(\xi)=\int_X f(x)\cdot
a_H(\xi^{-1}x)^\lambda \, dx\, .$$
We will call the map $\lambda\mapsto \widehat{f}_\lambda\in \cO(\Xi_+)$ the
{\it spherical holomorphic Fourier transform of $f$}.

\begin{lemma}\label{convergence} The following assertions
hold:
\begin{enumerate}
\item Let $U\subset A_+$ be a compact subset. The series
$$\sum_{\lambda\in \Lambda_{>0}} a_H(\xi^{-1})^\lambda \qquad (\xi\in\Xi_+)$$
converges uniformly on $GU\cdot \xi_0 \subset \Xi_+$.
\item For all $\xi\in \Xi_+$ one has
$$\sum_{\lambda\in \Lambda_{>0}} a_H(\xi^{-1})^\lambda =\cK(\xi)
\, .$$
\end{enumerate}
\end{lemma}

\begin{proof} Uniform convergence on $GU\cdot\xi_0$
is immediate from Corollary \ref{cor=1} and
Lemma \ref{lemma=pos}. Summing up the geometric series one obtains
\begin{eqnarray*} \sum_{\lambda\in \Lambda_{>0}} a_H(\xi^{-1})^\lambda &=&
\sum_{k_1=\ldots=k_{m-1}=0}^\infty\sum_{k_m=1}^\infty
a_H(\xi^{-1})^{k_1 \omega_1+\ldots+k_m \omega_m} \\
&=& \left({1\over
1-a_H(\xi^{-1})^{\omega_m}}-1\right)\cdot
\prod_{j=1}^{m-1}{1\over
1-a_H(\xi^{-1})^{\omega_j}}\\
&=& {1\over
a_H(\xi^{-1})^{-\omega_m}-1}\cdot
\prod_{j=1}^{m-1}{1\over
1-a_H(\xi^{-1})^{\omega_j}}\\
&=&\cK(\xi)
\end{eqnarray*}
\end{proof}

We conclude from Lemma \ref{convergence} that
the horospherical Cauchy transform of a function $f\in L^1(X)$ can be
decomposed as
$$\widehat{f}=\sum_{\lambda\in \Lambda_{>0}} \widehat{f}_\lambda $$
with the right hand side converging uniformly on compacta.

\begin{remark} We wish to point out that the
horospherical Cauchy kernel is a product of
geometrical and not  functional analytic reasoning.
We emphasize that in general not all parameters
$\lambda\in \Lambda_{>0}$ in the decomposition of the horospherical
Cauchy kernel correspond to unitarizable lowest weight modules
(see Remark \ref{rem=la} below for a more detailed
discussion).
\end {remark}

\subsection{Holomorphic Fourier transform on
lowest weight representations}
The objective of this subsection is to give a more
detailed discussion of the holomorphic Fourier transform for
functions $f\in L^2(X)$ which are contained in lowest weight module.
\par To begin with we collect  some material
on spherical unitary lowest weight representations. A reasonable source
might be the overview article \cite{ko}.

\par Let $(\pi_\lambda, {\mathcal H}_\lambda)$ be a non-trivial
 $H$-spherical unitary lowest weight
representation of $G$. As before we denote by $v_\lambda$ a normalized
lowest weight vector. Write $v_H$ for the unique $H$-fixed distribution
vector which satisfies $\langle v_\lambda, v_H\rangle =1$.

We record the fundamental
identity

\begin{equation}\label{eq=mc}
a_H(x)^\lambda=\langle \pi_\lambda(x)v_H, v_\lambda\rangle\qquad (x\in X)\, ,
\end{equation}
which allows us to link our geometric discussion in Section
\ref{section=hor-II} with representation theory.

\begin{remark} It follows from Corollary \ref{cor=D} that
$a_H^\lambda$ admits a holomorphic extension to the minimal
tube $D_-$. Traditionally this fact was explained via
(\ref{eq=mc}) in the context of holomorphic extension of
unitary lowest weight modules (see
\cite{n99}). We wish to point
out that  Corollary \ref{cor=D} asserts more,
namely that $a_H^\lambda|_{D_-}$ is bounded by $1$.
In addition Corollary \ref{cor=D} is more geometric, i.e.
not restricted to unitary parameters $\lambda$.
\end{remark}

Pairing the $G$-module of smooth vectors $\cH_\lambda^\infty $
with $v_H$ yields the $G$-equivariant embedding

\begin{equation}\label{eq=Xembed}
\iota: {\mathcal H}_\lambda^\infty \to C^\infty (X), \ \ v\mapsto
\left(x\mapsto
 \langle \pi_\lambda(x^{-1})v, v_H\rangle\right)\, .\end{equation}
We say that $\pi_\lambda$ is {\it $X$-square integrable} if
there exists a constant $d_s(\lambda)>0$, the
{\it spherical formal dimension} (cf. \cite{k}),
such that $\sqrt{d_s(\lambda)}\cdot \iota$
extends to an isometric map  $\cH_\lambda\to L^2(X)$.
\par $X$-square
integrable parameters $\lambda$ are characterized by the condition
\cite{oo91}
\begin{equation}\label{par2}
\la\lambda-\rho, \alpha\ra>0 \qquad \text {for all $\alpha\in \Delta_n^+$}\ .
\end{equation}
Here $\rho={1\over 2}\sum_{\alpha\in \Delta^+}
m_\alpha \alpha$ with $m_\alpha=\dim_\C \fg_\C^\alpha$.

\par Likewise we say $\pi_\lambda$ is {\it $X$-integrable} if
$\iota(\cH_\lambda^{\rm K-fin})\subset L^1(X)$. Integrability
is described by the inequality
\begin{equation}\label{par1}
\la\lambda-2\rho, \alpha\ra>0
\qquad \text {for all $\alpha\in \Delta_n^+$}\ .
\end{equation}

\par The set of parameters
$\lambda\in \Lambda_{>0}$ which satisfy condition (\ref{par1}), resp.
(\ref{par2}), shall be denoted by $\Lambda_1$, resp. $\Lambda_2$.
Note that $\Lambda_1\subset \Lambda_2$.

\begin{remark}\label{rem=la} We will
discuss the lattice $\Lambda_{>0}$ with regard to $\Lambda_1$ and
$\Lambda_2$. One recognizes a strong dependence
on the multiplicities  $m_\alpha$ which we will exemplify for
three basic cases below. Recall that elements $\lambda\in
\Lambda_{>0}$ are described by $\lambda=\sum_{i=1}^n \lambda_i
\omega_i$ with $\lambda_i\in \Z_{\geq 0}$ and $\lambda_m>0$.
In addition let us keep in mind that conditions
(\ref{par2}) and (\ref{par1}) are equivalent to
$\la\lambda-\rho, \alpha_m\ra>0$, resp.
$\la\lambda-2\rho, \alpha_m\ra>0$.

\par The equal rank case: In this situation one has $\ft=\fc$ and
$m_\alpha=1$ for all $\alpha$. Thus
$\rho={1\over 2}\sum_{i=1}^m \omega_i$ and therefore
$\lambda-\rho=\sum_{i=1}^m (\lambda_i-{1\over 2})\omega_i$.
In particular $\la \lambda-\rho,\alpha_m\ra =
\la \alpha_m, \alpha_m\ra (\lambda_m-{1\over 2})$ and thus
$\Lambda_{>0}\subset \Lambda_2$ as $\lambda_m\geq 1$ for elements
$\lambda\in \Lambda_{>0}$.

\par The group case: In this situation one has
$m_\alpha =2$ for all $\alpha$ and so
$\rho=\sum_{i=1}^n \omega_i$. Accordingly we obtain
$\la \lambda-\rho, \alpha_m\ra =
\la \alpha_m, \alpha_m\ra (\lambda_m-1)$.
It follows that $\Lambda_{>0}$ parameterizes
the holomorphic discrete series and their limits; in particular
$\Lambda_2 \subset \Lambda_{>0}$.

\par The rank one case:  Here one has $\Lambda_{>0}=\Z_{>0}\cdot \omega$
and $\rho={{m_\alpha}\over 2}\alpha$. Thus
$\Lambda_2=(\Z_{>0}+ \left[{m_\alpha\over 2}\right])\omega$ and
$\Lambda_2\subset \Lambda_{>0}$ with equality precisely for
$m_\alpha=1$, i.e. $\fg=\mathfrak{sl}(2,\R)$.
\end{remark}

\par For $\lambda\in \Lambda_2$ we set  $L^2(X)_\lambda=\iota(\cH_\lambda)$.

\begin{lemma}\label{lemma=Schur} Let $\lambda, \mu\in \Lambda_2$. Fix $v\in \cH_\lambda$
and define $f(x)=\la \pi_\lambda (x^{-1})v, v_H\ra\in L^2(X)_\lambda$.
Then for all $\xi\in \Xi_+$, the function
$$X\to \C,\ \  x\mapsto f(x)a_H(\xi^{-1}x)^\mu$$
is integrable and
\begin{equation}
\label{eq=Int} \int_X f(x)a_H(\xi^{-1}x)^\mu \, dx=
{\delta_{\lambda\mu}\over
d_s(\lambda)}
\la v, \pi_\lambda(\oline \xi)v_\lambda\ra\ .
\end{equation}
Here
$$\pi_\lambda(\oline \xi)v_\lambda=a^{-\lambda}
\pi_\lambda(g)v_\lambda\in \cH_\lambda \qquad\hbox{for
$\xi=ga\cdot \xi_0$, $g\in G$  and $a\in A_+$}\ .$$
\end{lemma}

\begin{proof} Fix $\xi\in \Xi_+$. Holomorphic extension of
(\ref{eq=mc}) yields
$$a_H(\xi^{-1} g)^\mu =\la \pi_\lambda(g)v_H, \pi_\lambda(\oline \xi)
v_\lambda\ra$$
for all $g\in G$ and $\xi\in \Xi_+$. It follows that
$x\mapsto a_H(\xi^{-1}x)^\mu$ is square integrable on $X$. Thus
$x\mapsto f(x)a_H(\xi^{-1}x)^\mu$ is integrable.
Finally we apply  Schur-orthogonality
(cf.\ \cite{k}, Prop. 3.2) and obtain

\begin{eqnarray*}
\int_{X}f(x) a_H(\xi^{-1}x)^\mu \, dx
 & =& \int_{X}\langle \pi_\lambda (x^{-1})v, v_H^\lambda\rangle
\langle \pi_\mu (x)v_H^\mu, \pi_\mu(\oline \xi)v_\mu\rangle
\, dx\\
& =& \int_{X}\langle \pi_\lambda (x^{-1})v, v_H^\lambda\rangle
\oline{\langle \pi_\mu (x^{-1})\pi_\mu(\oline {\xi} )v_\mu, v_H^\mu\rangle}
\, dx\\
& =& {\delta_{\mu\lambda}\over d_s(\lambda)}
\langle v, \pi_\lambda(\oline \xi)v_\lambda\rangle\, .
\end{eqnarray*}
\end{proof}

For  $\lambda\in \Lambda_1$ let us write
$L^1(X)_\lambda$ for the closure of
$\iota(\cH_\lambda^{\rm K-fin})$ in $L^1(X)$.
The next lemma
can be understood as an $L^1$-version of Schur-orthogonality
for the Cauchy-transform.

\begin{lemma}\label{lem1} Let $\lambda\in \Lambda_1$. Then
$$\widehat{f} =\widehat{f}_\lambda\qquad \hbox {\rm for all $f\in L^1(X)_\lambda$}\, .$$
\end{lemma}

\begin{proof}  Fix $f\in L^1(X)_\lambda$.
We have to show that $\widehat{f}_\mu=0$ for all
$\mu\in \Lambda_{>0}\backslash \{\lambda\}$. For
$\mu\in \Lambda_2$ this is a  consequence of
Lemma \ref{lemma=Schur}.
Therefore, we may assume that
$\mu\in \Lambda_{>0}\backslash \Lambda_2$. It means
that condition (\ref{par2}) is violated which
we will express as
\begin{equation} \label{parin}
\la \mu-\rho,\alpha_m\ra \leq 0\, .\end{equation}

\par We now show that
$$\widehat{f}_\mu(\xi)=\int_X f(x) a_H(\xi^{-1}x)^\mu\ dx =0 \qquad
\text{for all $\xi\in \Xi_+$}\, .$$
Equation above has boundary values on $G/M\subset \partial \Xi_+$
and it will be sufficient to prove that
$$\widehat{f}_\mu(gM)=\int_X f(x) a_H(g^{-1}x)^\mu\ dx =0 \qquad
\text{for all $g\in G$}\, .$$
We compute
\begin{eqnarray*}
\widehat{f}_\mu(gM)&=&\int_X f(x) a_H(g^{-1}x)^\mu\ dx\\
&=&\int_X f(gx) a_H(x)^\mu\ dx\\
&=&\int_T \int_X f(tgx) a_H(tx)^\mu\ dx\, dt\\
&=&\int_X \left(\int_T t^\mu f(tgx)\,  dt\right)
 a_H(x)^\mu\ dx\, .
\end{eqnarray*}
To arrive at a contradiction, suppose that
$\int_T t^\mu f(tgx)\, dt\neq 0$. This
can only happen if $\mu$ belongs to the $T$-weight spectrum of $\pi_\lambda$.
Now the $T$-weights of $\pi_\lambda$ are contained
in $\lambda+\Z_{\geq 0}[\Delta^+]$. Thus $\mu=\lambda+\gamma$ for some
$\gamma\in \Delta^+$.
But then
$$\la \mu-\rho, \alpha_m\ra =\la \lambda-\rho, \alpha_m\ra +
\la \gamma, \alpha_m\ra\, .$$
Observe that both summands on the right hand side
are positive, the desired contradiction to (\ref{parin})
\end{proof}

\begin{remark}\label{rem0} (Analytic continuation) Let $\lambda\in \Lambda_1$ and $f\in L^1(X)_\lambda
\cap L^2(X)_\lambda$. Write $f(x)=\la \pi_\lambda(x^{-1})v, v_H\ra $
for some $v\in \cH_\lambda$.
Then Lemma \ref{lemma=Schur} and Lemma \ref{lem1} imply
\begin{equation}\label{eq=44}
\widehat{f}(\xi) ={1\over d_s(\lambda)}
\la v, \pi_\lambda(\oline \xi)v_\lambda\ra\, .
\end{equation}
Clearly, the right hand side makes sense for all
$v\in \cH_\lambda$ and all $X$-square integrable parameters
$\lambda\in \Lambda_2$. We now explain how passing
to parameters $\lambda\in \Lambda_2$ in (\ref{eq=44})
has a natural explanation in terms of
analytic continuation.
For that let $\tilde G$ denote the universal cover of $G$. Write
$\tilde \Lambda_1$,
$\tilde \Lambda_2$ for the sets of $\tilde G$-integral parameters
which satisfy (\ref{par1}), resp. (\ref{par2}).
Clearly $\Lambda_{1,2}\subset \tilde\Lambda_{1,2}$. The effect
of passing to the universal cover is that the parameter
spaces involved become continuous in the central variable,
i.e. there exists constants $0<c_2<c_1$ such that
$$\tilde \Lambda_1|_{\R Z_0}=]c_1,\infty[\cdot (\alpha_m|_{\R Z_0})
\quad\hbox{and}\quad\tilde \Lambda_2|_{\R Z_0}=
]c_2,\infty[\cdot (\alpha_m|_{\R Z_0})\ .$$
By the concrete formula for $d_s(\lambda)$ from \cite{k}, Th. 4.15,
we know that $\lambda\mapsto d_s(\lambda)$ is a meromorphic function
on $\fa_\C^*$ which is  positive on $\tilde \Lambda_2$.
Now familiar techniques show that
that the assignment   $\tilde\Lambda_2\ni\lambda\mapsto{1\over d_s(\lambda)}
\la v, \pi_\lambda(\oline \xi)v_\lambda\ra\in \C $ becomes
analytic in the central variable (the Shapovalov
form is polynomial in $\lambda$ and in
 \cite{k99} it is explained how to make consistent analytic
choices for $v$ and $v_\lambda$
in dependence of the central coordinate of $\lambda$).
\end{remark}

Motivated by Remark \ref{rem0}
we define the horospherical Cauchy transform
for functions $f\in L^2(X)_\lambda$, $\lambda\in \Lambda_2$ by
$$\widehat{f}=\widehat{f}_\lambda\, .$$

\subsection{Hyperfunctions and generalized matrix coefficients}
In order to discuss the horospherical Cauchy transform and its inverse
in a more comprehensive  way
we need some results on the analytic continuation of
generalized matrix coefficients of lowest weight representations. Proofs
of the facts cited below can be found in \cite{KNO}.

\par Let $(\pi, \cH)$ be a unitary
lowest weight representation of $G$.  Write
$\cH^\omega$ and $\cH^{-\omega}$ for the associated
$G$-modules of analytic, resp. hyperfunction vectors.
The nature of the $T$-spectrum of $\pi$ shows that
$\pi|_{T}$ extends holomorphically to $\cT_-$. Moreover
the so obtained self adjoint operators
$\pi(a)$, $a\in A_-$, are of trace class  and strongly mollifying, i.e
\begin{equation}\label{eq=moll}
\pi(a)\cH^{-\omega}\subset \cH^\omega \qquad (a\in A_-)\, .\end{equation}

\par Assume that $\pi$ is $H$-spherical and denote
by $v_H$ the (up to scalar) unique $H$-fixed distribution vector.
Let $v\in \cH^{-\omega}$ be a hyperfunction vector.
We wish to interpret
the generalized matrix coefficient $f(x)=\la\pi(x^{-1})v, v_H\ra$
as a generalized function on $X=G/H$.
It follows essentially from (\ref{eq=moll})
that the prescription
\begin{equation} \label{eq=he}
\tilde f(ga\cdot x_0)=\la \pi(g^{-1})\pi(a^{-1})v, v_H\ra \qquad
\hbox{for $g\in G$ and $a\in A_+$}\ .\end{equation}
defines a holomorphic function on $D_+=GA_+\cdot x_0$.
The minimal tube $D_+$ has $X$ as edge
and this allows us to interpret $f$ as the boundary
value of $\tilde f$. Henceforth we will identify
$f$ with the holomorphic function $\tilde f$.

\par Suppose that $\cH \subset L^2(X)$, i.e. $\pi=\pi_\lambda$
with $\lambda\in \Lambda_2$ and $\cH=L^2(X)_\lambda$.
We now show how the horospherical Cauchy transform restricted
to $L^2(X)_\lambda$ can be extended to $L^2(X)_\lambda^{-\omega}\subset
\cO(D_+)$. In other words  for $\xi\in \Xi_+$
we wish to give
meaning to
\begin{eqnarray*}
\widehat{f}(\xi)&=&\int_X f(x) a_H(\xi^{-1}x)^\lambda\ dx\\
&=& \int_X \la \pi(x^{-1})v, v_H\ra a_H(\xi^{-1}x)^\lambda\ dx
\end{eqnarray*}
as a holomorphic function on $\Xi_+$.
Express $\xi$ as $\xi=ga\cdot \xi_0$ with
$g\in G$ and $a\in A_+$. By the usual holomorphic change
of variables  one obtains that

\begin{eqnarray*}
\widehat{f}(\xi)&=&\int_X f(x) a_H(a^{-1}g^{-1}x)^\lambda \ dx \\
&=& \int_X \la \pi(x^{-1})\pi(g^{-1})\pi(a^{-1})v, v_H\ra a_H(x)^\lambda
\ d(x)\ .
\end{eqnarray*}
Now the last expression is well defined by (\ref{eq=moll}).
Of course one has
\begin{equation}\label{eq=lala}
\widehat{f}(\xi)={1\over d_s(\lambda)} \la \pi(a^{-1}) \pi(g^{-1})
v, v_\lambda\ra  \end{equation}
by the same argument as in Lemma \ref{lemma=Schur}.
Thus we have shown that the horospherical Cauchy transform on $L^2(X)_\lambda$
extends to a $G$-equivariant continuous map
$$L^2(X)_\lambda^{-\omega}\to \cO(\Xi_+)_\lambda\,. $$

\par We conclude this section with a conjecture
related to the holomorphic intertwining of $\cO(D_+)$ and
$\cO(\Xi_+)$. It can be seen as a holomorphic analogue of
Helgason's conjecture (actually a theorem ny \cite{K-}).
\par In order to state the conjecture some new terminology is
needed.
Let us call a holomorphic function
$f$ on $\Xi_+$ {\it bounded away from the boundary} if
its restriction to $g\cT_+a$ is bounded for all
choices of $g\in G$ and $a\in A_+$. We denote by
$\cO_{\rm b.a.b.}(\Xi_+)$ the space of all holomorphic function
on $\Xi_+$ which are bounded away from the boundary. Note that
$\cO_{\rm b.a.b.}(\Xi_+)$ is a closed $G$-subspace of the
Fr\'echet space $\cO(\Xi_+)$.

\begin{conjecture}\label{con=1} Let ${\Bbb D}(X)$ be the algebra
of $G$-invariant differential operators on $X$. Naturally
we can view ${\Bbb D}(X)$ as holomorphic differential
operators on $X_\C$. Write $\cO(D_+)_\lambda$ for the
common holomorphic ${\Bbb D}(X)$-eigenfunctions on $D_+$ with infinitesimal
character $\lambda-\rho$. Let $\lambda\in \Lambda_2$.
We conjecture
\begin{equation}
\cO_{\rm b.a.b.}(D_+)_\lambda=L^2(X)_\lambda^{-\omega}\ . \end{equation}
Notice that the inclusion $''\supset''$ is clear by
(\ref{eq=moll}).
\par We have already remarked that $\cO(\Xi_+)_\lambda\simeq
\cH_\lambda^{-\omega}$ \cite{KNO}. Hence our conjectured equality
means that the horospherical Cauchy transform induces
an intertwining isomorphism $\cO_{\rm b.a.b.}(D_+)_\lambda\to \cO(\Xi_+)_\lambda$.
\par It is also interesting problem to formulate (and prove)
the conjecture for other parameters.
\end{conjecture}

\begin{remark} We illustrate Conjecture \ref{con=1} for the one-sheeted
hyperboloid $X= {\rm Sl}(2,\R)/ {\rm SO} (1,1)$.
Fix $\lambda\in 2\Z_{>0}=\Lambda_2$ and denote by
$\cH_\lambda$, resp. $\cH_{-\lambda}$,  the lowest (resp. highest)
weight module of $G={\rm Sl}(2,\R)$ with lowest (resp. highest) weight
$\lambda$, resp. $-\lambda$. Denote by $V_{\lambda-2}$ the
finite dimensional $G$-module of highest weight $\lambda-2$.
Write $C^\infty(X)_{\lambda}$ for the
${\Bbb D}(X)$-eigenspace   with eigenvalue $\lambda-1$.
Then
$$C^\infty(X)_\lambda \simeq
\cH_\lambda^\infty \oplus \cH_{-\lambda}^\infty \oplus V_{\lambda-2}$$
Now, the functions of $\cH_{\pm \lambda}^\infty $ extend holomorphically
to $D_\pm$ but not beyond,  while the functions
of $V_{\lambda-2}$ extend holomorphically to all of $X_\C$.
One deduces that $\cO(D_+)_\lambda=\cH_\lambda^{-\omega}\oplus
V_{\lambda-2}$. Finally, the holomorphic
functions in $V_{\lambda-2}$ grow exponentially at infinity and hence
are not bounded away from the boundary.
Thus $\cO_{\rm b.a.b.}(D_+)_\lambda =L^2(X)_\lambda^{-\omega}$ as conjectured.
\end{remark}

\subsection{Inversion of the horospherical Cauchy transform} To begin with
we first have to explain certain facts on
incidence geometry between the Shilov boundary $X$ of $D_+$ and the
boundary piece $G/M$ of $\Xi_+$.
\par We keep in mind that we realized
$G/M$ in the boundary of $\Xi_+$ by
$G/M\simeq  G\cdot \xi_0\subset \partial \Xi_+$.

\par  Recall the orbits $S(z)\subset \Xi$ from
\ref{eq=S}. For a point $x\in X$ we
define the real form of $S(x)$ by

$$S_\R(x)=S(x)\cap G/M\ .$$
In view of the incidence relation (\ref{eq=indi}),
one has
$$S_\R(x)=\{ \xi\in G/M\mid \xi\in S(x)\}=
\{ \xi\in G/M\mid x\in E(\xi)\}\,  .$$

\begin{lemma}\label{lemma=fiber} Let $x=g\cdot x_0\in X$, $g\in G$. Then
$$S_\R(x)=gH\cdot \xi_0\simeq H/M\ .$$
\end{lemma}

\begin{proof} First notice that for $x=g\cdot x_0$ with $g\in G$ one has
$S_\R(x)= g\cdot S_\R(x_0)$. Hence it suffices to show that
$$S_\R (x_0)=H\cdot \xi_0\simeq H/M\ .$$
Let $\xi\in S_\R(x_0)$ and write $\xi=y\cdot \xi_0$ for some
$y\in G$. We have to show that $y\in H$ and that
$y$ is uniquely determined modulo $M$. First observe that
$\xi\in S(x_0)$ means $yN_\C \subset H_\C N_\C$
and so
$y\in H_\C N_\C \cap G$.
Now $H_\C N_\C \cap G=H$ implies $y\in H$. Finally, uniqueness modulo $M$
is immediate from Lemma \ref{l-12one}.
\end{proof}

It is possible to view the boundary orbits  $S_\R(x)$ as
certain limits.
For $z=ga\cdot x_0\in D_+$ we define
$$S_\R(z)=gaH\cdot \xi_0\simeq H/M\, .$$
We note that $S_\R(z)\subset \Xi_+$ by
(\ref{newxi}). Furthermore there is the
obvious limit relation
$$\lim_{a\to {\bf 1}\atop a\in A_+} S_\R(ga\cdot\xi)=S_\R(g\cdot x_0)\, .$$

\par Write $d_z(\xi)$ for the measure on $S_\R(z)$ which is
induced from a Haar measure $d(hM)$ on $H/M$ via the identification
$S_\R(z)\simeq H/M$. Define the space of {\it fiber integrable}
holomorphic functions on $\Xi_+$ by

$$\cO_{\rm f.i.}(\Xi_+)=\{ \phi\in\cO(\Xi_+)\mid  D_+ \ni z\to
\int_{S_\R(z)} |\phi(\xi)|\,  d_z(\xi) \quad \hbox{is locally bounded}\}\ .$$

For a function $\phi\in \cO_{\rm f.i.}(\Xi_+)$ we define its
{\it inverse horospherical transform} $\phi^\vee\in \cO(D_+)$
by
$$\phi^\vee(z)=\int_{S_\R(z)} \phi(\xi) \, d_z(\xi) \qquad
(z\in D_+)\, .$$
We note that
$$\cO_{\rm f.i.}(\Xi_+)\to \cO(D_+), \ \ \phi\mapsto \phi^\vee$$
is a $G$-equivariant continuous map.

\par  Finally,  we define a subset $\Lambda_c\subset \Lambda_2$ of large
parameters by

$$\Lambda_c=\{ \lambda\in \Lambda_2\mid (\forall \alpha\in \Delta_n^+)
\ (\lambda -\rho)(\check\alpha)> 2 - m_\alpha\}\, . $$

The inversion formula for the horospherical Cauchy transform is based on the following key result.

\begin{lemma} Let $\lambda\in \Lambda_c$. Let
$f\in L^2(X)_\lambda^{-\omega}\subset\cO(D_+)$.
Then $\widehat{f}\in \cO_{\rm f.i}(\Xi_+)$
and
\begin{equation}\label{eq=inversion}
f(z)=d(\lambda)\cdot  \int_{S_\R(z)} \widehat{f}(\xi) \ d_z(\xi)
\qquad (z\in D_+)\, .\end{equation}
In other words, $f=d(\lambda)\cdot (\widehat{f})^\vee$.
\end{lemma}

\begin{proof} Let $f(ga\cdot \xi_0)=
\langle \pi_\lambda(a^{-1})\pi_\lambda(g^{-1})v, v_H\rangle$
for some $v\in {\mathcal H}_\lambda^{-\omega}$.
Then by (\ref{eq=lala})

$$\widehat{f}(\xi) ={1\over d_s(\lambda)}
\langle v, \pi_\lambda (\oline{\xi})v_\lambda\rangle\, .$$
As $\lambda\in \Lambda_c$, \cite{k}, Th. 2.16 and Th. 3.6, imply that
$$\int_{H/M}\pi_\lambda(h)v_\lambda \, d(hM)= {d_s(\lambda)\over d(\lambda)}\cdot v_H$$
with the left hand side understood as convergent
$\cH_\lambda^{-\omega}$-valued integral.

Thus with $z=ga\cdot x_0$ one obtains that

\begin{eqnarray*} \int_{S_\R(z)}\widehat{f}(\xi)\, d_z(\xi) & =& \int_{H/M}
{1\over d_s(\lambda)}
\langle \pi_\lambda(a^{-1})\pi_\lambda (g^{-1})v, \pi_\lambda(h)v_\lambda\rangle
\, d(hM)\\
&=& {1\over d(\lambda)}\cdot  f(z)\, ,
\end{eqnarray*}
completing the proof of the lemma.
\end{proof}

\begin{remark} If $\lambda\in \Lambda_2
\backslash \Lambda_c$
and $0\neq f\in L^2(X)_\lambda^{-\omega}$, then
the integral $\int_{S_\R(z)} \widehat{f}(\xi)\, d_z(\xi)$ does not converge.
However, using the results from \cite{k} it can be shown
that the identity (\ref{eq=inversion}) can be
analytically continued (cf. Remark \ref{rem0})
to all $\lambda\in \Lambda_2$.
Henceforth we understand  (\ref{eq=inversion})
as an identity valid for all $\lambda\in \Lambda_2$.
\end{remark}
The formal dimension $d(\lambda)$ is a polynomial in $\lambda$,
explicitly given by
\cite{hc55}
\begin{equation}\label{***}d(\lambda)=
 c \cdot
\prod_{\alpha\in \Sigma^+} \langle  \lambda-\rho(\fc), \alpha\rangle
\end{equation}
with $c\in \R$ a constant depending on the normalization
of measures.

\par The right action of $T$ on $\Xi_+$ induces an identification
of ${\mathcal U} (\ft_\C)$ with
$G$-invariant differential operators on $\Xi_+$.
As usual we identify ${\mathcal U}(\ft_\C)$ with polynomial functions on
$\ft_\C$. In this
way $d(\lambda)$ corresponds to a $G$-invariant
differential operator $\cL$ on $\Xi_+$ which acts
along the fibers of $\Xi_+\to F_+$ and has constant
coefficients in logarithmic coordinates.
In particular,
$$\cL \phi = d(\lambda) \cdot \phi  \qquad (\phi\in \cO(\Xi_+)_\lambda)\ .$$
Combining this fact with equation (\ref{eq=inversion}) we obtain the
main result of this paper.

\begin{theorem} Let
$f\in \sum_{\lambda\in \Lambda_2} L^2(X)_\lambda^{-\omega}\subset \cO(D_+)$.
Then
$$ f = (\cL \widehat{f})^\vee\ .$$
\end{theorem}

\section{The example of the hyperboloid of one sheet}
\noindent
This section is devoted to the discussion of the
case $G={\rm Sl}(2,\R)$ and $H={\rm SO}(1,1)$. Notice that
$G/H\simeq {\rm SO}_e(2,1)/{\rm SO}_e(1,1)$. For what follows it
is inconsequential to assume that $G={\rm SO}_e(2,1)$ and $H={\rm
SO}_e(1,1)$
although the universal complexification of
$G={\rm SO}_e(2,1)$ is not simply connected.

\par The map
$$G/H\to \R^3, \ \ gH\mapsto g\cdot\left[\begin{matrix}1\\ 0\\ 0\end{matrix}\right]$$
identifies $X=G/H$ with the one sheeted hyperboloid
$$X=\{ x=(x_1, x_2, x_3)^T\in \R^3\mid x_1^2 + x_2^2 -x_3^2 =1\}\ .$$
The base point $x_0$ becomes $(1,0,0)^T$.
Let us define a complex bilinear pairing on $\C^3$ by
$$\la z,  w\ra = z_1w_1+ z_2w_2 -z_3w_3 \qquad \text{for}\quad
z=\left[\begin{matrix}z_1\\ z_2\\ z_3\end{matrix}\right], w=\left[\begin{matrix}w_1\\ w_2\\ w_3\end{matrix}\right]\in \C^3\, .$$
If we set $\Delta(z)=\la z, z\ra $ for $z\in \C^3$, then
$X=\{ x\in \R^3\mid \Delta(x)=1\}$.
Further one has $G_\C={\rm SO}(2,1; \C)\simeq {\rm SO}(3,\C)$
and $H_\C={\rm SO(1,1; \C)}\simeq {\rm SO}(2,\C)$. Clearly
$$X_\C=G_\C/ H_\C =\{ z\in \C^3\mid \Delta(z)=1\}\ .$$
\par Our choice of $T$ will be
$$T=K=\left\{ \left[\begin{matrix} \cos \theta & \sin \theta & 0\\  -\sin\theta & \cos \theta & 0\\
0 & 0 & 1\end{matrix}\right]\mid \theta\in \R\right\}\ .$$
In particular $\fa=\R U_0$ where
$$U_0=\left[\begin{matrix} 0  & i & 0\\  -i  & 0 & 0\\
0 & 0 & 0 \end{matrix}\right]$$
and $\Delta=\Delta_n=\{ \alpha, -\alpha\}$ with $\alpha(U_0)=1$. If we
demand $\alpha$ to be the positive root, then

$$N_\C =\left\{ \left[ \begin{matrix} 1-{z^2\over 2} & i{z^2\over 2} & iz\\
i{z^2\over 2}  & 1+{z^2\over 2} & z\\
iz & z & 1\end{matrix}\right]\mid z\in \C\right\}\ .$$
The homogeneous space $G_\C/M_\C N_\C $ naturally
identifies with the isotropic vectors $\Xi=\{\zeta\in
\C^3\backslash\{0\} \mid \Delta(\zeta)=0\}$
via the $G_\C$-equivariant map
$$G_\C/ M_\C N_\C\to \Xi,\ \  gM_\C N_\C\mapsto g\cdot\zeta_0 \qquad \text{where}
\quad \zeta_0=\left[\begin{matrix} 1 \\  -i\\  0\end{matrix}\right].$$
The correspondence between elements
of $\zeta\in \Xi$ and horospheres on $X_\C$ is
explicitly given by
$$\zeta \leftrightarrow E(\zeta)=\{ z\in X_\C\mid \la z,  \zeta\ra =1\}\ .$$
Elements $\zeta\in \Xi$ can be expressed as $\zeta=\xi+i \eta$ with
$\xi, \eta\in \R^3\backslash\{0\}$ subject to
$$\Delta(\xi)=\Delta(\eta) \qquad \text{and}\qquad \la \xi, \eta\ra =0\, .$$
A simple computation yields

$$\Xi_+=\{ \zeta=\xi+i\eta\in \Xi\mid \Delta(\xi)=\Delta(\eta)>1\}\, $$
and
$$D_+=\{ z=x+iy\in X_\C \mid \Delta(x)> 1\}\, .$$

Next we compute the kernel function.

\begin{lemma}\label{lem=ah} For all $z\in X_\C$ and $\zeta\in \Xi$
one has
$$a_H(\zeta^{-1} z)^{-\alpha}= \la z,\zeta\ra\ . $$
\end{lemma}

\begin{proof} We first show that
\begin{equation}\label{eq=33}
a_H(g)^{-\alpha}= \la g\cdot x_0, \zeta_0\ra\qquad (g\in G_\C)\ .
\end{equation}
Observe that both sides are holomorphic functions on $G_\C$
which are left $N_\C$ and right $H_\C$-invariant.
Thus it is enough to test with
elements $a\in A
_\C$. Then $a_H(a)^{-\alpha}=
a^{-\alpha}$.
On the other hand
for
$a= \left[\begin{matrix} \cos \theta & \sin \theta & 0\\  -\sin\theta & \cos \theta & 0\\
0 & 0 & 1\end{matrix}\right]$ with $\theta\in \C$ we specifically obtain

$$\la a\cdot x_0, \zeta_0\ra =\la \left[\begin{matrix} \cos\theta\\ -\sin\theta\\ 0\end{matrix}
\right],  \left[\begin{matrix} 1 \\ -i\\ 0\end{matrix}
\right]\ra= \cos\theta +i \sin\theta =a^{-\alpha}\ .$$
This proves (\ref{eq=33}).
\par It is now easy to prove the asserted statement
of the lemma. For that write $\zeta=g\cdot \zeta_0$ and $z=y\cdot x_0$
for $g,y\in G_\C$.
Then  (\ref{eq=33}) implies

$$a_H(\zeta^{-1}z)^{-\alpha}=a_H(g^{-1}y)^{-\alpha}=
\la g^{-1}y\cdot x_0, \zeta_0\ra  =\la y\cdot x_0,  g\cdot\zeta_0\ra=
\la z, \zeta\ra\,.$$
\end{proof}

We observe that $\Lambda_{>0}=\Lambda_2=\Z_{>0}\cdot\alpha$. Hence
Lemma \ref{lem=ah}
implies that the horospherical Cauchy kernel is
$$\cK(\zeta)={1\over a_H(\zeta^{-1})^{-\alpha}-1}= {1\over \la \zeta,
 x_0\ra  -1}\qquad
(\zeta\in \Xi_+)\ .$$

The horospherical Cauchy transform for $f\in L^1(X)$ is given by

$$\widehat{f}(\zeta)=\int_X {f(x)\over \la \zeta,  x\ra  -1}\  dx \qquad
(\zeta\in \Xi_+)$$
with $dx$ the invariant measure  on the hyperboloid $X$.
Finally, we will discuss inversion. Let the inner product
on $\fa$ be normalized such that $\langle \alpha, \alpha\rangle =1$
and identify $\R$ with $\fa^*$ by means of the bijection
$\R\ni \lambda\mapsto \lambda\alpha\in \fa^*$. Then $\Lambda_{>0}=\Z_{>0}$
and $d(\lambda)=\lambda-{1\over 2}$. An easy calculation gives

$$\cL=\sum_{j=1}^3 \zeta_j {\partial\over \partial \zeta_j} \ -\ {1\over 2}
\ .$$
For $f\in \sum_{\lambda>0} L^2 (X)_\lambda^{-\omega}\subset \cO(D_+)$ the inversion formula
reads
$$f(z)=\int_{-\infty}^\infty (\cL f)\left(\begin{matrix}z_1 -i {z_2\over r
}\cosh t -i {z_1z_3\over r}\sinh t\\
 z_2 +i{z_1\over r}\cosh t
-i{z_2z_3\over r }\sinh t\\
z_3 -i r\sinh t\end{matrix} \right) \ dt,  $$
where $r=\sqrt{z_1^2+z_2^2}$.

\end{document}